%% file: article_arXiv.tex
\documentclass[12pt]{article}

\usepackage{amsfonts,amsmath,amssymb,epsfig,epsf,psfrag,bbm}
\usepackage[latin1]{inputenc}
\usepackage{graphicx}
\usepackage{color}
\usepackage{xcolor}
\usepackage{url}

\RequirePackage{natbib}

\input{commands.tex}

\begin{document}
\title{Degree-based goodness-of-fit tests
 for heterogeneous random graph models :
independent and exchangeable cases}
%
%

\author{
  Sarah Ouadah$^{1, *}$, St\'ephane Robin$^1$, Pierre Latouche$^2$ \\
  \footnotesize{($^1$) UMR MIA-Paris, INRA, AgroParisTech, Universit\'e Paris-Saclay, 75005 Paris, France} \\
  \footnotesize{($^2$) Universit\'e de Paris, MAP5, CNRS, 75006 Paris,
    France}  
  }


\maketitle

\begin{abstract}
The degrees are a classical and relevant way to study the topology of a network. 
They can be used to assess the goodness-of-fit for a given random graph model.
In this paper we introduce goodness-of-fit tests for two classes of models.
First, we consider the case of independent graph models such as the heterogeneous Erd\"os-R\'enyi model in which the edges have different connection probabilities. 
Second, we consider a generic model for exchangeable random graphs called the $W$-graph. The stochastic block model and the expected degree distribution model fall within this framework.
We prove the asymptotic normality of the degree mean square under these independent and exchangeable models and derive formal tests. 
We study the power of the proposed tests and we prove the asymptotic normality under specific sparsity regimes.
The tests are illustrated on real networks from social sciences and ecology, and their performances are assessed via a simulation study.
\end{abstract}

\textbf{Keywords:} degree variance; goodness-of-fit ; graphon; random graphs; $W$-graph.

\maketitle

\section{Introduction}
Interaction networks are used in many fields such as biology,
sociology, ecology, economics or energy to describe the interactions
existing between a set of individuals or entities. Formally, an
interaction network can be viewed as a graph, the nodes of
which being the individuals, and an edge between two nodes
being present if these two individuals
interact. Characterizing the general organization
of such a network, namely its topology, can help in understanding the behavior of the system as a whole.

In the last decades, the distribution of the degrees (i.e. the number
of connections of each node) has appeared
as a simple and relevant way to study the topology
of a network, see \cite{Sni81} and \cite{BaA99}. The degree distribution can also be used to infer complex graph models such as in \cite{BCL11}. From
a more descriptive view-point, a very imbalanced distribution may
reveal a network whose edges highly concentrate around few nodes, whereas a multi-modal
distribution may reveal the existence of clusters of nodes as observed by
\cite{CDR12}. However, in
  practice, assessing the significance of such
patterns remains an open problem. \\ 
The variance of the
degrees has been considered since the earliest statistical studies of
  networks, for instance by \cite{Sni81}. The first idea was simply to compare its
empirical value to the expected one under a null
random graph model, typically the Erd\"os-R\'enyi ($ER$) model introduced by \cite{ER59}, where each degree has a
binomial distribution.
Because the
$ER$ model is rarely a reasonable model to be tested, we define a
generalized version of the degree variance statistic, which we name
the degree mean square statistic. This statistic generalizes the degree
variance in the sense that it measures the discrepancy between the
observed degrees and their expected values under several heterogeneous models we define hereafter.\\
For a given random graph under a specific model $M_0$, the degree mean square statistic is defined by
\[W_{\theta_0}= \frac1n \sum_i (D_i - \Esp_{\theta_0} D_i)^2,\]
where $D_i$ stands for the degree of node $i$ and $\Esp_{\theta_0} D_i$ for its expected value under a given model with parameter $\theta_0$. More specifically, in the following we will consider independent models parametrized with a probability matrix $\pg$ and exchangeable models parametrized with a function $\Phi$. Although these models not only differ in terms of parameter, for the sake of clarity, the corresponding quantities will be simply indexed with $\pg$ and $\Phi$, respectively.
We propose 
goodness-of-fit tests for several random graph models, by showing the asymptotic normality of this statistic $W_{\theta_0}$ under null hypothesis and their alternatives. In addition, because large networks are often sparse, we study under which sparsity regime the asymptotic
distributions derived before still hold.

The notations and the main models considered are the following. We consider an undirected graph $\Gcal =
(\{1, \dots n\}, \mathcal{E})$ with no self loop, that is
  the connection of a node to itself, and denote $Y$ the
corresponding $n \times n$ adjacency matrix. Thus, the entry
$\Yij$ of $Y$ is 1 if $(i, j) \in \mathcal{E}$, and 0
otherwise. Because $\Gcal$ is undirected with no self loop, we have $\Yij = Y_{ji},\forall i\neq j$ and $Y_{ii} = 0$, for all $i$'s. We further define $D_i$ the degree of
node $i$ by $D_i = \sum_{j \neq i} Y_{ij}$. 
In terms of random graph
models, we consider two cases: the independent case and the exchangeable one.
In the independent case, $ER(p)$ refers to the Erd\"os-R\'enyi model, according to which
all edges $(\Yij)$ are independent Bernoulli variables with same
probability $p$ to exist. $HER(\pg)$ stands for the heterogeneous
Erd\"os-R\'enyi model where edges are independent with respective
probability $\pij$ to exist. The $n\times n$ matrix $\pg$ has entries
$\pij$, it is symmetric with null diagonal. 
In the exchangeable case,
we consider a generic model for exchangeable random graphs called the $W$-graph introduced in \cite{LoS06} and \cite{DiJ08}. It is based on a {\sl graphon} function $\Phi: [0, 1]^2 \mapsto [0, 1]$ and denoted by $\EGphi$. An unobserved coordinate $U_i \sim \mathcal{U}[0, 1]$ is associated with each node $i (1\leq i\leq n)$ and edges are drawn independently conditional the $U_i$'s as
$
Y_{ij}|U_i, U_j \sim \Bcal[\Phi(U_i, U_j)].
$
The stochastic block model (SBM, introduced by \cite{Hol79} and further studied by \cite{NoS01}, and the expected degree distribution (EDD) model, defined by \cite{ChL02}, fall within this framework. 

Goodness-of-fit tests of the models we consider have received little attention until recently.
\cite{CFV15} propose a goodness-of-fit test for the $HER(\pg)$ model and\cite{Mau17} for the $\EGphi$ model, both when independent and identically distributed (i.i.d.) copies of the graph are available. \cite{Lei16} and \cite{BiSa16} derived goodness-of-fit tests for the number of communities in stochastic block models by showing the asymptotic behavior of the largest singular value of a residual adjacency matrix. Their respective null models are $ER(p)$ in \cite{BiSa16} and an SBM with $K$ communities in \cite{Lei16}. \cite{YHA14} proposed a test statistic for the goodness-of-fit of a given graphon function and used a Monte-Carlo sampling to approximate its null distribution. More recently, \cite{Gao17} proved the asymptotic normality of subgraph counts to test the $ER(p)$ model against an SBM.

The paper is organized as follows.
Section 2 is devoted to independent graph models and Section 3 to the the exchangeable ones. The performances of the proposed tests are assessed via a simulation study in Section 4. More specifically, the asymptotic distribution of the degree mean square statistic under models $HER(\pg)$ and $\EGphi$ is derived Sections \ref{Sec:MeanSqDeg} and \ref{Sec:MeanSqDeg:Graphon}, respectively. The asymptotic normality under some
specific sparsity regimes is studied in Sections \ref{Sec:Sparse} and \ref{Sec:Sparse:Graphon}. In Section \ref{Sec:TestDMS}, we establish a test for the null hypothesis stating that $\Gcal$ arises from $\HERpgo$ and give its power. The last part of this section is devoted to the illustration of the HER goodness-of-fit test on some examples. In the same manner, Section \ref{Sec:TestDMS:Graphon} deals with the EG model and its extensions, meaning the SBM and EDD model. 
\section{Independent random graph models \label{Sec:Indep}}

We consider the heterogeneous Erd\"os-R\'enyi model $HER(\pg)$, in which the edges are
independent and have different respective probabilities to exist:
$Y_{ij} \sim \Bcal[\pij]$.

The asymptotic framework is the following.
\begin{assumption}\label{Assumption}
In the non-sparse setting, we consider an infinite matrix ${\bf P}$, the elements of
which are all in the interval $[c,1-c]$ for some arbitrarily small constant $c\in(0,1/2)$. For the HER model, then we build a sequence
of matrices $\pg^n$ made of the first $n$ rows and columns of ${\bf
  P}$. Finally, we consider a sequence of independent graphs $\Gcal^n
= (\{1, \dots n\}, \mathcal{E}^n)$, with increasing size $n$ and respective
probability matrices $\pg^n$.  
\end{assumption}
All quantities computed
on $\Gcal^n$ should therefore be indexed by $n$ as well. But 
for the sake of clarity,
we will drop the index $n$ in $\pg^n$ in the rest of the paper.

\subsection{Asymptotic normality}   \label{Sec:MeanSqDeg}
We consider a goodness-of-fit test for the  $\HERpgo$ model. For a given random graph with a matrix $\pgo$ of connection probabilities, we consider the following degree mean square statistic:
\[\Wo = \frac1n \sum_i (D_i - \muo_i)^2,\]
where $D_i=\sumji\Yij$ and $\muo_i$ stand for the expected degree of node $i$ under $\HERpgo$, namely $\muo_i = \sumji \poij$.

We establish the asymptotic normality of $\Wo$ under model $HER(\pg)$. The proof relies on projections of $\Wo$ on suitable spaces and the Lindeberg-L\'evy Theorem (see e.g. Theorem 7.2, p.42 in \cite{Bil68}) which is recalled below. We derive all projections involved in the Hoeffding decomposition (see, e.g., Chapter 11 in \cite{vdV98}) to easily calculate the moments of $\Wo$. As for the asymptotic normality, we decompose $\Wo$ into the sum of its H\'ajek projection (see, e.g., Chapter 11 in \cite{vdV98}) to which we apply the Lindeberg-L\'evy Theorem, and a negligible term. A similar strategy has already been used for graph studies, for instance in \cite{Blo05} to prove the asymptotic normality of the variance degree under model $ER(p)$ and in \cite{NoW88} to prove the one of subgraph counts in random graphs.

\begin{theorem}[Lindeberg-L\'evy Theorem in \cite{Bil68}] \label{Thm:LindebergLevy}
{\it  Let  $(\Xnu)_{1\leq u\leq  k_n}$ be a triangular array of independent random variables with means $0$ and finite variances $(\sigmanu^2)_{1\leq u\leq k_n}$. Let $B_n^2=\sum_{u=1}^{k_n}\sigmanu^2$.
If the Lindeberg condition
\begin{eqnarray}\label{lind}
{A_n^2(\epsilon)}/{B_n^2} \to 0, \quad \mbox{as } n\to\infty, \quad \mbox{for each }\epsilon>0,
\quad \text{where} \quad 
A_n^2(\epsilon) = \sum_{u=1}^{k_n} \int_{\{|x_{nu}|>\epsilon B_n\}} x_{nu}^2 dP
\end{eqnarray}
is satisfied then
\begin{eqnarray*}
\frac1{B_n} \sum_{u=1}^{k_n}\Xnu  \overset{\mathcal{D}}{\longrightarrow} \mathcal{N}(0, 1).
\end{eqnarray*}}
\end{theorem}

\begin{remark}\label{rem1} Let consider the case of binary random variables $\Xnu$ with mean 0. More specifically, set $\Xnu=\anu\Znu$, $\anu\in\mathbb{R}$, where $\Znu$ are centered Bernoulli variables, that is to say $\Znu$ takes value $1-p_{nu}$ with probability $\pnu$ and value $-\pnu$ with probability $1-\pnu$. Because $|\Xnu| \leq \anu$, the realization of the event $|\Xnu| \geq \epsilon B_n$ in the definition of $A_n^2(\epsilon)$ in \eqref{lind} is controlled by $|\anu|\geq \epsilon B_n$. Therefore, all $\Xnu$ for which $|\anu| < \epsilon B_n$ do not contribute to $A_n^2(\epsilon)$. If this holds for all $\Xnu$, then the Lindeberg condition is directly satisfied. If not, only the $\Xnu$ for which it does not hold have to be considered in the calculation of $A_n^2(\epsilon)$ and, because $|\Znu| \leq 1$, their contribution is upper-bounded by their variance $\sigmanu^2 = \anu^2 \pnu(1-\pnu)$. In the forthcoming theorems proofs, we will verify the Lindeberg condition using this observation. \end{remark}

\begin{theorem} \label{Thm:W0AsympNormHER}
 Under model $HER(\pg)$ and Assumption \ref{Assumption}, the statistic $\Wo$ is asymptotically normal:
 $$
 (\Wo - \Esp_{\bf{p}} \Wo)/\Sd_{\bf{p}} \Wo \overset{D}{\longrightarrow} \Ncal(0, 1),
 $$
where $\Sd$ denotes the standard deviation and
\begin{equation}
\Esp_{\bf{p}} \Wo = \frac{2}{n} \left( \sum_{1 \leq i < j \leq n} (\sigmaij + \dcij) + \sum_{1 \leq i < j < k \leq n} (\dij\dik + \dij\djk + \dik\djk) \right),\label{Thm:W0AsympNormHER:Esp}
\end{equation} 
where $\sigmaij=\pij(1-\pij)$ and $\dij = \pij - \poij$. Moreover
 \begin{eqnarray}\label{Thm:W0AsympNormHER:Var}
  \Var_{\bf{p}} \Wo & =& \frac{4}{n^2} \Big(\sum_{1\leq i<j\leq n} \sigmaij (1 - 2 \pij+\Di + \Dj )^2\nonumber\\
&&+ \sum_{1\leq i<j<k\leq n} \left(\sigmaij\sigmaik + \sigmaij\sigmajk + \sigmaik\sigmajk \right) \Big),
 \end{eqnarray}
with $\Delta_i = \sumji \dij$.
\end{theorem}
\proofbegin
Let us begin with the calculation of $\Wo$ moments. We
first observe that, 
\begin{eqnarray*} 
n \Wo
& = & \sum_i (D_i - \mu_i + \mu_i - \mu^0_i)^2 
\; = \; \sum_i \left(\sum_{j\neq i} \Ytij + \dij\right)^2 \\
& = & 2 \sum_{1 \leq i < j \leq n} (\Ytij + \dij)^2 \\
& & + 
2 \sum_{1\leq i<j<k\leq n} (\Ytij + \dij)(\Ytik + \dik) +  (\Ytij + \dij)(\Ytjk + \djk) +  (\Ytik + \dik)(\Ytjk + \djk),
\end{eqnarray*} 
where $\Ytij=\Yij-\pij$ and $\mu_i = \sumji \pij$.
Then, we write the Hoeffding decomposition of $\Wo$:
\begin{eqnarray}\label{hoef2}
\Wo = 
P_{\emptyset} \Wo 
+ \sum_{1\leq i<j\leq n} P_{\{ij\}}\Wo
+ \sum_{1\leq i<j<k\leq n} \left(P_{\{ij,ik\}} \Wo + P_{\{ij,jk\}} \Wo + P_{\{ik,kj\}} \Wo \right),
\end{eqnarray}
where
\begin{eqnarray*}
P_{\emptyset}\Wo&=&\Esp_{\bf{p}}\Wo,\\
P_{\{ij\}}\Wo&=&\Esp_{\bf{p}}(\Wo|\Yij)-\Esp_{\bf{p}}\Wo,\\
P_{\{ij,ik\}}\Wo&=&\Esp_{\bf{p}}(\Wo|\Yij,\Yik)-\Esp_{\bf{p}}(\Wo|\Yij)-\Esp_{\bf{p}}(\Wo|\Yik)+\Esp_{\bf{p}}\Wo.
\end{eqnarray*}
Combining the definitions above with the expression \eqref{hoef2} of $\Wo$, we obtain that,
\begin{eqnarray}
P_\emptyset \Wo &=& \frac2n \sum_{1 \leq i < j \leq n} (\sigmaij + \dcij) + \frac2n \sum_{1 \leq i < j < k \leq n} \dij\dik + \dij\djk + \dik\djk, \nonumber\\
P_{\{ij\}} \Wo &=& \frac2n \Ytij \left(1+\Di + \Dj \right)-\sigmaij= \frac2n \Ytij \left(1-2\pij + (\Di + \Dj) \right),\label{pij2}\\
P_{\{ij, ik\}} \Wo &=& \frac2n \Ytij\Ytik.\label{pijk2}
\end{eqnarray}
Observe now that,
\begin{eqnarray*}
n\Esp_{\bf{p}}\Wo&=& 2\sum_{1 \leq i < j \leq n} (\sigmaij + \dcij) + 2 \sum_{1 \leq i < j < k \leq n} \dij\dik + \dij\djk + \dik\djk.
\end{eqnarray*}
Because the $\Ytij$ are independent with zero mean, the projections are all orthogonal with each other, which gives
\begin{eqnarray*}
n^2 \Var_{\bf{p}}\Wo & = & n^2 \sum_{1 \leq i < j \leq n}\Var_{\bf{p}}(P_{\{ij\}} \Wo)\\
&&+ n^2 \sum_{1 \leq i < j < k \leq n}\left(\Var_{\bf{p}} (P_{\{ij, ik\}} \Wo)+\Var_{\bf{p}}(P_{\{ij, jk\}} \Wo)+\Var_{\bf{p}}(P_{\{ik, jk\}} \Wo)\right)\\
&=& 4\sum_{1 \leq i < j \leq n}\sigmaij (1 - 2 \pij + \Di + \Dj)^2 +4\sum_{1\leq i<j<k\leq n} \left(\sigmaij\sigmaik + \sigmaij\sigmajk + \sigmaik\sigmajk \right).
\end{eqnarray*} 

We now turn to the asymptotic normality of $\Wo$.
Let decompose $\Wo$ as follows.
\begin{eqnarray*}
\Wo-\Esp_{\bf{p}}\Wo=\Wo^*-\Esp_{\bf{p}}\Wo+\Wo-\Wo^*,
\end{eqnarray*}
where $\Wo^*=P_{\emptyset}\Wo+\sum_{1\leq i<j\leq n}P_{\{ij\}}\Wo$ is the H\' ajek projection of $\Wo$, which corresponds to the first two terms of the Hoeffding's decomposition. We will show that $\Wo^*-\Esp_{\bf{p}}\Wo$ is asymptotically normal and that $\Wo-\Wo^*$ is a negligible term.\\
Let consider $\Wo^*-\Esp_{\bf{p}}\Wo=\sum_{1\leq i<j\leq n}P_{\{ij\}}\Wo$ and apply Theorem \ref{Thm:LindebergLevy} to the projections $P_{\{ij\}}\Wo$ which stand for the $\Xnu$. We first observe that these projections are each proportional to the $\Ytij$ which are all independent centered Bernoulli variables. We may now use Remark \ref{rem1}. We denote the $\anu$ by $a_{n\{ij\}}$, the explicit expression of which is given in \eqref{pij2}. We observe that, under Assumption \ref{Assumption}, $a_{n\{ij\}}=\Theta(1)$ and $B^2_n=\Var_{\bf{p}}\left(\Wo^*-\Esp_{\bf{p}}\Wo\right)= \Theta(n^2)$. It implies that the Lindeberg condition is fulfilled because, for any $\epsilon$, each $\anu$ becomes smaller than $\epsilon B_n$ when $n$ goes to infinity. Now by considering \eqref{hoef2} the Hoeffding decomposition of $\Wo$, we see that $$\Wo-\Wo^*= \sum_{1\leq i<j<k\leq n} \left(P_{\{ij,ik\}} \Wo + P_{\{ij,jk\}} \Wo + P_{\{ik,kj\}} \Wo \right).$$ Then we observe that $a_{n\{ij,ik\}}$ given in \eqref{pijk2} is $\Theta(n^{-1})$ and therefore that $\Var_{\bf{p}}\left(\Wo-\Wo^*\right)=\Theta(n)$. We conclude to the asymptotic normality of $\Wo$ by combining the one of $\Wo^*-\Esp{\bf{p}}\Wo$ and the fact that $\Var_{\bf{p}}\left(\Wo-\Wo^*\right)/\Var_{\bf{p}}\Wo^*\to 0$ as $n\to\infty$.
\proofend 

\paragraph{Plug-in version of the test.}
In many situations, $\pgo$ is actually unknown and one needs to resort to an estimate $\pgoh$. There is no hope to get a precise estimate when $n$ increases if no restriction is imposed to $\pgo$. When a vector of covariates $\xij \in \mathbb{R}^d$ is available for each pair of nodes, one natural way to impose such a restriction is to assume that $\poij$ has a logistic form, that is $\text{logit}(\poij) = x_{ij}^\intercal \beta$, where $\beta$ is the vector of regression coefficients and $\text{logit}(u) = \log(u/(1-u))$. 
A plug-in version of the proposed test can be obtained by fitting the logistic model to the observed edges to get an estimate $\widehat{\beta}$, which provides us with $\pgoh$, which in turn provides us with a plug-in version $\Who$ of the test statistic. \\
The simulation study presented in Section \ref{Sec:Simus} shows that $\Who$ behaves well for large graphs. 
A possible strategy to understand the asymptotic behavior of $\Who$ would be to control the difference between $\Wo$ and $\Who$. Indeed, denoting $\muho_i = \sum_{j \neq i} \phoij$ and $\Delta_i = \muho_i-\muo_i$, $\Who$ can be decomposed as
\begin{equation} \label{Eq:Who}
\Who 
:= \frac1n \sum_i \left(D_i - \muho_i\right)^2
= \Wo - \frac2n \sum_i \left(D_i - \muo_i\right) \Delta_i + \frac1n \sum_i \Delta_i^2.
\end{equation}
If the $\phoij$ result from a parametric estimation based on the $O(n^2)$ edges, we expect the estimation error $|\poij - \phoij|$ to be $O_P(n^{-1})$, which makes the last term of \eqref{Eq:Who} negligible. Still, the joint dependence structure of the $D_i$ and $\Delta_i$ is quite intricate, which makes the control of the second term of \eqref{Eq:Who} not straightforward. 
In Section \ref{Sec:TestDMS}, we present a specific case where we prove the asymptotically normality of the plug-in version of the test.

\subsubsection*{Degree variance test} \label{Sec:VarDeg}

We consider the following statistic
which is the empirical degree variance for the test of $H_0 = ER$ versus $H_1 = HER(\pg)$.
\[\Vh = \frac1n \sum_i \left(D_i - \Db\right)^2,\] 
where $\Db=(1/n)\sum_{j}D_j$. The variance of the
degrees has been naturally considered earlier in statistical studies of networks. \cite{Hag03a} derives the exact moments of
the degree variance and suggests to use a Gamma distribution as in
\cite{Hag00}. \cite{Sni81} also gives the first two moments of the
degree variance, but conditionally to the total number of edges. To
our knowledge the first and only proof of the asymptotic normality of
the degree variance under the ER model is given in a technical report
from \cite{Blo05}. Here, we establish the asymptotic normality of $\Vh$ under model $HER(\pg)$ and obtain the ER version as a consequence.

\begin{corollary} \label{Thm:VhAsympNormHER}
 Under model $HER(\pg)$ and Assumption \ref{Assumption}, 
the degree variance is asymptotically normal:
 $$
 \left(\Vh-\Esp_{\bf{p}} \Vh\right) / \Sd_{\bf{p}}\Vh \overset{D}{\longrightarrow} \Ncal(0, 1),
 $$
with
\begin{eqnarray*}
\Esp_{\bf{p}}\Vh 
& = 
&\frac{2(n-2)}{n^2}\sum_{1\leq i<j\leq n} \pij 
+ \frac{2(n-4)}{n^2}\sum_{1\leq i<j<k\leq n} \left\{\pij\pik+\pij\pjk+\pik\pjk\right\}\\
&&-\frac{8}{n^2}\sum_{1\leq i<j<k<l\leq n}\left\{\pij \pkl+\pik\pjl+\pil\pjk\right\},\end{eqnarray*}
and
\begin{eqnarray*}
\Var_{\bf{p}}\Vh&=&\frac{1}{4n^4}\sum_{1\leq i<j\leq n} \sigmaij\left(4(n-2)+4(n-4)\sum_{k\notin(i,j)}(p_{i,k}+p_{j,k})-16\sum_{k<l\notin(i,j)}\pkl\right)^2\\
&&+\frac{1}{n^4}\sum_{1\leq i<j<k\leq n} 4(n-4)^2\left\{\sigmaij\sigma^2_{ik}+\sigmaij\sigmajk+\sigma^2_{ik}\sigmajk\right\}\nonumber\\
&&+\frac{1}{n^4}\sum_{1\leq i<j<k<l\leq n}
64\left\{\sigmaij\sigmakl+\sigma^2_{ik}\sigmajl+\sigmail\sigmajk\right\}.
\end{eqnarray*}
\end{corollary}

The proof follows the line of this of Theorem \ref{Thm:W0AsympNormHER} and is given in Appendix \ref{Proof:Thm:VhAsympNormHER}.\\
Note that the asymptotic normality of the degree variance under model $ER(p)$ is a straightforward application of Corollary \ref{Thm:VhAsympNormHER} to the case where all $\pij$ are equal to $p$. We have,
$$
\left(\Vh - \Esp_{p} \Vh\right)  / \Sd_{p}\Vh \overset{D}{\longrightarrow} \Ncal(0, 1),
$$
where $\Esp_{p}\Vh = n^{-1} {(n-1) (n-2) p q}$ and $\Var_{p}\Vh =  n^{-3} {2 (n-1) (n-2)^2} p q \left(1+(n-6) pq \right)$, as given in \cite{Hag00}.

\subsection{Test and power} \label{Sec:TestDMS}
We now study the test of $H_0 = \HERpgo$ versus $H_1 = \HERpg$. The next Corollaries provide the null distribution of the test statistic $\Wo$ and the power of the associate test.

\begin{corollary} \label{Cor:W0AsympNormHER0}
 Under model $\HERpgo$ and Assumption \ref{Assumption}, the statistic $\Wo$ is asymptotically normal with moments:
  \begin{align}
  \Esp_{\bf{p^0}} \Wo & = \frac2n \sum_{1 \leq i < j \leq n} \sigmaoij, \label{Cor:W0AsympNormHER0:Esp}\\
  \Var_{\bf{p^0}} \Wo & = \frac1{n^2} \left(4\sum_{1\leq i<j\leq n} \sigmaoij (1 - 2 \poij)^2
+ \sum_{1\leq i<j<k\leq n} \left(\sigmaoij\sigmaoik + \sigmaoij\sigmaojk + \sigmaoik\sigmaojk \right) \right),\label{Cor:W0AsympNormHER0:Var}
 \end{align}
where $\sigmaoij=\poij(1-\poij)$.
\end{corollary}

This is a direct consequence of Theorem \ref{Thm:W0AsympNormHER} in the special case of the $\HERpgo$ model for which all $\dij$'s are zero ($\dij = \pij - \poij$).

 A formal test with asymptotic level $\alpha$ can be constructed based on Corollary \ref{Cor:W0AsympNormHER0}, which rejects $H_0$ as soon as $\Wo$ exceeds $\Esp_{\bf{p^0}}\Wo + t_{1-\alpha} \Sd_\HERpgo \Wo$, where $t_{1-\alpha}$ stands for the $1-\alpha$ quantile of the standard Gaussian distribution. The power of this test is given by the following Corollary.

\begin{corollary}
The asymptotic power of the test for $H_0 = HER(\pg^0)$ versus $H_1 = HER(\pg)$ is 
\begin{align}
 \pi(\pg) = 1 - \Phi\left(\left(\Esp_{\bf{p^0}} \Wo + t_{1-\alpha} \Sd_{\bf{p^0}} \Wo - \Esp_{\bf{p}} \Wo \right) \left/ \Sd_{\bf{p}} W_\pgo \right.\right),\label{Cor:HER:Power}
\end{align}
where $\Phi$ stands for the cumulative distribution function [cdf] of the standard normal distribution and $t_{1-\alpha} = \Phi^{-1}(1-\alpha)$.
\end{corollary}

The following corollary gives a sufficient condition on the departure between $\pg$ and $\pg^0$ to ensure that the proposed test is asymptotically powerful.
\begin{corollary} \label{Cor:W0AsympNormHER0:Pwer}
For probability matrices $\pg^0$ and $\pg$, define
\begin{eqnarray*}
\Delta_n(\pg^0, \pg) &:=& \Esp_{\pg} \Wo - \Esp_{\pg0} \Wo\\
&=&\frac{2}{n} \left( \sum_{1 \leq i < j \leq n} \dcij + \sum_{1 \leq i < j < k \leq n} (\dij\dik + \dij\djk + \dik\djk) +\sum_{1 \leq i < j \leq n} (\sigmaij-\sigmaoij)\right).
\end{eqnarray*}
If $\Delta_n(\pg^0, \pg) = \Theta(n^\alpha)$ is positive and $\alpha > 1/2$, then under Assumption \ref{Assumption}, the test $H_0 = HER(\pg^0)$ versus $H_1 = HER(\pg)$ is asymptotically powerful.
\end{corollary}
\proofbegin
It is sufficient to prove that the argument of the cdf $\Phi$ in \eqref{Cor:HER:Power} tends to minus infinity as $n$ increases. From \eqref{Thm:W0AsympNormHER:Var} and \eqref{Cor:W0AsympNormHER0:Var}, we have that under Assumption \ref{Assumption}, $\Sd_{\pg} \Wo = \Theta(n^{1/2})$ and $\Sd_{\pg0} W_\pgo = \Theta(n^{1/2})$. As a consequence,when $\Delta_n(\pg^0, \pg)>0$ and $\alpha > 1/2$, the negative argument of $\Phi$ in \eqref{Cor:HER:Power} goes to infinity at rate $n^{\alpha - 1/2}$, which concludes the proof.
\proofend

Note that, when $\Delta_n(\pg^0,\pg)<0$ the same corollary holds for the test which rejects $H_0$ as soon as $\Wo<\Esp_{\bf{p^0}}\Wo + t_\alpha \Sd_\HERpgo \Wo$, where $t_{\alpha}$ stands for the $\alpha$ quantile of the standard Gaussian distribution.

\subsubsection*{Degree variance test}\label{Sec:TestDV}
We now consider the use of the statistic $\Vh$ for the test of $H_0 = ER$ versus $H_1 = HER(\pg)$. 
Because the  probability is unknown in practice, we consider the following test statistic using a plug-in version of the moments, namely
\[ 
\left(\Vh - \Esp_{\widehat{p}} \Vh\right)  / \Sd_{\widehat{p}} \Vh,
\]
where $\ph=[n(n-1)]^{-1} \sumij\Yij$. \\
The asymptotic power $\pi(\pg)= \Pr_\pg\{\Vh > t_\alpha\}$ of the considered test, with nominal level $\alpha>0$, is
\[\pi(\pg) = 1 - \Phi\left(\left(\Esp_{\bar{p}} \Vh + t_{1-\alpha} \Sd_{\bar{p}}\Vh - \Esp_{\bf{p}} \Vh\right) \left/ \Sd_{\bf{p}}\Vh \right.\right),\]
where $\pb=[n(n-1)]^{-1}\sumij\pij$.
This results from the asymptotic normality of $(\Vh - \Esp_{\bar{p}} \Vh)  / \Sd_{\bar{p}}\Vh$ under the $HER(\pg)$ model. Actually, the asymptotic distribution of the test based on $(\Vh - \Esp_{\bar{p}} \Vh)  / \Sd_{\bar{p}}\Vh$ is the same as the one of the test based on the statistic $(\Vh - \Esp_{\widehat{p}} \Vh)  / \Sd_\ERph\Vh$ (see Lemma \ref{Lem:VphEquivVpb} in Appendix \ref{Lemma:PowerER}), and we have shown that under model ER, $(\Vh - \Esp_{\widehat{p}} \Vh)  / \Sd_\ERph\Vh$ is asymptotically normal (see Lemma \ref{Cor:VhTestER} in Appendix \ref{Lemma:PowerER}).
\begin{remark} 
The $ER(p)$ model corresponds to $HER(\pgo)$ where the matrix $\pgo$ has all entries equal to $p$. In this case, the test statistic $\Wo$ can be viewed as the theoretical version of the empirical variance statistic $\Vh$ studied in Section \ref{Sec:VarDeg} as
$$
\Wo = \frac1n \sum_i \left(D_i - (n-1)p\right)^2.
$$
Because as $\ph$ is an average over $\Theta(n^2)$ edges, we have that $(\ph - p)^2 = \Theta_P(n^{-2})$ so $\Wo - \Vh = (n-1)^2 (\ph - p)^2 = \Theta_P(1)$. Combined with arguments similar to these of Corollary \ref{Cor:VhTestER} and Lemma \ref{Lem:VphEquivVpb}, this implies that, under the ER model, the tests based on $\Vh$ and $\Wo$ are asymptotically equivalent.
\end{remark}

\subsubsection*{Illustration}\label{Sec:Illustration}
We illustrate the use of the proposed test on the following series of networks.
\begin{description}
 \item[{\sl Karate network}:]  it describes the friendships between a subset of $n=34$ members of a karate club at a university in the  US, observed from 1970 to 1972 and was originally  studied by \cite{zachary1977}. 
The network is made of
four known groups characterized by a node
  qualitative descriptor.

 \item[{\sl Ecological networks}:] this consists in two ecological networks
   first   introduced  in \cite{VPD08}   and   further  studied   in
   \cite{MRV10}. Each of these networks describe the interaction between a series of $n = 51$ trees and $n = 154$ fungi, respectively. In the tree network, two trees interact if they share at least one common fungal parasite. As for the fungal network, two fungi are linked if they are hosted by at least one common tree species.
Three
{quantitative edge descriptors}
are available  characterizing the  genetic, geographic,  and taxonomic
distances between the tree species.

 \item[{\sl Political blogs network}:]  this consists in a set of $n = 196$ French political blogs studied in \cite{LBA11a}. Two blogs are connected if one contains an hyperlink to the other.

Each node is associated with a political party from the left
wing to  the right  wing and the  status of the  writer is  also given
(political analyst or not).


 \item[{\sl CKM}:] this data set was created  by \cite{Burt1987} from the data originally
collected   by   \cite{Coleman1966}.   The   network   we   considered
characterizes  the  friendship  relationships among $n=219$ physicians,  each
physician being asked to name  three friends. 

The physicians were also
asked to answer to a series of questions regarding their profession, corresponding to node covariates.
Note that we imputed
the missing values in the data set using the missMDA R package of \cite{husson2016}.
 \item[{\sl Faux Dixon High network}:]
this  network  characterizes  the (directed)  friendship  between  $n=248$
students.  It results  from  a simulation  based  upon an  exponential
random graph model fit, see \cite{handcock2008}, to data from one school community from the AdHealth Study, Wave I of \cite{Resnick97}.

 Node
covariates  are provided,  namely the  grade,  sex, and  race of  each
student. 

  \item[{\sl AdHealth 67}:] 
this data  set is  related to the  Faux Dixon  network described
  previously. However,  it was constructed  from the original  data of
  the  AdHealth  study,  and  not simulated  from  any  random  graph
  model.   The  AdHealth   study   was   conducted  using   in-school
  questionnaires, from 1994  to 1995. Students were  asked to designate
  their friends and  to answer to a series of  questions. Results were
  collected in schools from 84 communities. In our study, we considered a network
  associated to school community 67 which characterizes the undirected
  friendship relationships between $n=530$ students. 

Nodes covariates are the same as the one of the Faux Dixon
  network.
 \end{description}

For some networks, only node descriptors $x_i$ and $x_j$ are available and building an {\sl edge} descriptor $x_{ij}$ from {\sl node} descriptors  is not straightforward as depicted in \cite{HGH08}. In these examples, the node descriptors are all qualitative. For each category of each node  descriptor, we build binary edge descriptors indicating if both node belong to the same category, or if at least one on the two belong to it.
The precise definition of the edge covariates for each dataset is explained in \cite{gof2}.


 We fist applied the degree variance test to each of these networks to check if their topology is similar to the one of an ER network. As  expected, their topology  are far
too heterogeneous to fit an $ER(p)$ model, and the  null hypothesis is
rejected for each one of them.\\
The question is then to know if the available covariates on edges are sufficient to explain the heterogeneity of the network, at least in terms of degrees. To address this question, for each network separately, we fitted a logistic regression model 
$\text{logit}(\poij) = x_{ij}^\intercal \beta$, which 
provided us with an estimate $\pgoh$ of the connection probability matrix $\pgo$. We then applied the degree mean square test to check if the considered covariates are sufficient to explain the heterogeneity of the network.

\begin{table}[!ht]
  \begin{center}
     \caption{Degree mean square HER test. 
     TestStat $= ({\Who-\Esp_{\bf{\widehat{p}^0}}})/{\Sd_{\bf{\widehat{p}^0}}}$. \label{Tab2:TreeFungi}}
\begin{tabular}{lccccccc}
 Network & $n$ & mean($\ph^0_{ij}$) & st-dev($\ph^0_{ij}$) & $\Who$ & $\Esp_{\bf{\widehat{p}^0}}\Who$ & $\Sd_{\bf{\widehat{p}^0}}\Who$ & TestStat \\ \hline

	Karate & 34&  0.135 & 0.149 & 3.84 & 3.22 & 0.88 & 0.71\\ 
    Trees & 51 & 0.553 & 0.2 & 140.23 & 10.66 &  2.11 & 61.55 \\ 
    Fungis & 154 & 0.226 & 0.021 & 592.12 & 26.82 &  3.06 & 184.55 \\
    Blogs & 196&  0.075 & 0.112 & 84.82 & 11.05 & 1.2 & 61.5\\  
CKM & 219&  0.015 & 0.035& 3.16 & 3 & 0.32 & 0.5\\
	Faux Dixon & 248&  0.02 & 0.037 & 11.34 & 4.41 & 0.43 & 16.05\\  
	AdHealth  & 530 &  0.007 & 0.008 & 8.77 & 3.43 & 0.24 & 22.27
    \end{tabular}
  \end{center}
\end{table}

The results given in Table \ref{Tab2:TreeFungi} show the ability of the proposed test to detect a departure from the degrees predicted by the covariates. 
Indeed, the null hypothesis is rejected for all networks except for the CKM and Karate networks. 
As for the ecological networks, these results are consistent with these from \cite{MRV10}, who detected a residual heterogeneity in the valued versions of these networks after correction for these covariates.


\subsection{Case of sparse graphs} \label{Sec:Sparse}
We discuss the validity of Theorem \ref{Thm:W0AsympNormHER} when considering sparse graphs. 
Sparsity can be defined in two ways. Either each connection probability vanishes as $n$ grows, or the fraction of non-zero connection probabilities decreases as $n$ grows. The following Proposition deals with a combination of both scenarios. 
\begin{proposition}\label{Prop:W0Sparse}
Consider the $HER(\pg)$ model, when $\pij=\psij n^{-a}$, $a>0$, $\psij$ following Assumption \ref{Assumption} and a fraction $1 - n^{-b}$, $b\geq 0$, of $\pij$'s is set to zero. The $\poij$'s satisfy the same assumptions. Then, provided that $a+b<2$, the statistic $\Wo$ is asymptotically normal.
\end{proposition}
\proofbegin 
We will show that $\Wo^*-\Esp_{\bf{p}}\Wo$ is asymptotically normal then that $\Wo-\Wo^*$ is a negligible term. The projections $P_{\{ij\}}\Wo$ involved in $\Wo^*-\Esp_{\bf{p}}\Wo$ still stand for the $\Xnu$ and $a_{n\{ij\}}$ expressed in \eqref{pij2} stand for $a_{nu}$ (notation of Remark \ref{rem1}). Since under Assumption \ref{Assumption} $\Di=\Theta(n^{1-a-b})$, we see that $a_{n\{ij\}}=\Theta(n^{-(a+b)})$ if $a+b<1$ and $\Theta(n^{-1})$ if $a+b>1$. Therefore, we have 
$\Var_{\bf{p}} P_{\{ij\}} V=\Theta\left(n^{-3a-2b}\right)$ if $a+b<1$ and $\Theta\left(n^{-a-2}\right)$ if $a+b>1$. Combining this with the number of non-zero terms which equals $\Theta(n^{2-b})$, we get that $B_n^2 =\Theta\left(n^{2-3(a+b)}\right)$ if $a+b<1$ and $\Theta\left(n^{-(a+b)}\right)$ if $a+b>1$.
Comparing $A_n^2(\epsilon)$ with $B_n^2$, we see that the Lindeberg condition is fulfilled for $a+b < 2$.\\
Now we consider $\Wo-\Wo^*$ as the sum of the projections $P_{\{ij,ik\}}\Wo$. The $a_{n\{ij,ik\}}$ given in \eqref{pijk2} equal $\Theta(n^{-1})$, thus $\Var_{\bf{p}} P_{\{ij,ik\}} \Wo=\Theta\left(n^{-2a-2}\right)$. Since the number of non-zero terms in the sum is $\Theta(n^{3-2b})$, we have therefore $\Var_{\bf{p}}\left(\Wo-\Wo^*\right)=\Theta(n^{1-2(a+b)})$.\\
We conclude to the asymptotic normality of $\Wo$ by combining the one of $\Wo^*-\Esp_{\bf{p}}\Wo$ under condition $a+b < 2$ and the fact that $\Var_{\bf{p}}\left(\Wo-\Wo^*\right)/\Var_{\bf{p}}\Wo^*\to 0$ as $n\to\infty$ under the same condition.
\proofend
\begin{remark}
The condition $a+b<2$ ensures that, although the density of the graph goes to zero, the number of edges still goes to infinity as $n$ grows.
\end{remark}

We now extend Corollary \ref{Thm:VhAsympNormHER} for the degree variance to sparse graphs, considering a setting similar to this of Proposition \ref{Prop:W0Sparse}.

\begin{corollary}\label{Prop:VhAsympNormHER}
Consider the $HER(\pg)$ model, with exactly the same conditions as in Proposition \ref{Prop:W0Sparse}. Then, provided that $a+b < 2$, the $V$ statistic is asymptotically normal.
\end{corollary}

The proof follows the line of this of Proposition \ref{Prop:W0Sparse} and is given in Appendix \ref{Proof:Prop:VhAsympNormHER}.

\section{Exchangeable random graph models}
We consider $\EGphi$ a generic model for exchangeable random graphs based on a {\sl graphon} function $\Phi: [0, 1]^2 \mapsto [0, 1]$ and commonly called the $W$-graph introduced in \cite{LoS06} and \cite{DiJ08}. Under $\EGphi$, a coordinate $U_i \sim \mathcal{U}[0, 1]$ is associated with each node $i (1 \leq i \leq n)$ and edges are drawn independently conditional the $U_i$'s as
\[
Y_{ij}|U_i, U_j \sim \Bcal[\Phi(U_i, U_j)].
\]
Many statistical models  such as the expected degree-corrected SBM, see \cite{dasgupta2004spectral,karrer11}, and the random Rash model, see \cite{rasch1960probabilistic}, fall into this framework. In this paper, we focus on the stochastic block model (SBM) and the expected degree distribution (EDD) model. 

\begin{description}
 \item[SBM.] The SBM introduced in \cite{Hol79} and \cite{NoS01} consists in a mixture model for random graph as pointed out by \cite{DPR08}, in which a discrete variable $Z_i \in\{1, \dots K\}$ is associated with each node and edges are drawn conditionally as $Y_{ij}|Z_i, Z_j \sim \Bcal[\pi_{Z_i,Z_j}]$, where $[\pi_{k\ell}]_{k, \ell}$ stands for the so-called connectivity matrix. Indeed, SBM corresponds to a $W$-graph with block-wise constant graphon function, see \cite{LaR15}.
 \item[EDD.] The EDD model is an exchangeable version of the expected degree sequence model studied in \cite{ChL02} and of the configuration model from \cite{New03}. Under these two models, the degree of each node is fixed which makes them non exchangeable. Under the EDD, an expected degree $K_i$ (not necessarily integer) is first drawn independently and identically for each node from some distribution $G$ and the edges are drawn independently conditional on the $K_i$ as $Y_{ij}|K_i, K_j \sim \Bcal[K_i K_j / \kappa]$, so $\Esp(D_i | K_i) \propto K_i$. EDD corresponds to a $W$-graph with product-form graphon function: $\Phi(u, v) = g(u) g(v)$, taking $g(u) = G^{-1}(u)/\sqrt{\kappa}$. \cite{YoS07} consider a specific case of this model.
\end{description}

\subsection{Asymptotic normality}   \label{Sec:MeanSqDeg:Graphon}
We propose a goodness-of-fit test for the $W$-graph model. For a given graphon $\phio$, we consider the following degree mean square statistic.
\[
\Wphio = \frac1n\sum_i (D_i -(n-1)\phioo)^2,\]
where $\phioo$ stands for the marginal probability for any given edge to exist, namely $\phioo=\int\int\phio(u,v)\dd u \dd v$.
We establish the asymptotic normality of $\Wphio$ under model $\EGphi$. The proof relies on a central limit theorem for acyclic patterns from \cite{BCL11}, which is recalled hereafter.

\begin{figure}[!ht] 
\begin{center}
\includegraphics[scale=0.85]{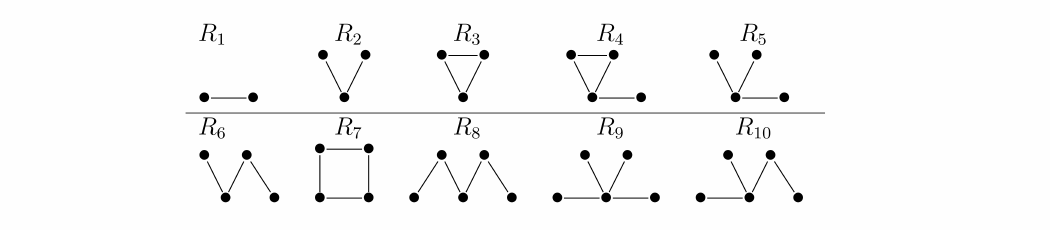}
 \caption{Definition of the patterns $R_1$ to $R_{10}$ involved in the calculation of the moment of the $W$ statistics.\label{Fig:Motifs}}
\end{center}
\end{figure}

Let us consider a fixed pattern $R$ (i.e. a given graph as displayed in Figure \ref{Fig:Motifs}) with $m$ nodes and set of edges $\Ecal_R$. Let us consider a random graph $\Gcal_R$ with $m$ nodes 
generated by $\EGphi$.  We define $P(R)$ and its empirical version $\hat{P}(R)$ computed on a graph $\Gcal$ with $n$ nodes as follows.
\begin{eqnarray}\label{phat}
P(R)=\mathbb{P}\left(\Gcal_R=R\right),\quad\mbox{ and }\quad
\hat{P}(R) = \dbinom{n}{m}^{-1} {N(R)}^{-1} \sum_{\Gcal_S\subset\Gcal} \mathbb{1}\left(\Gcal_S\sim R\right),
\end{eqnarray}
where $\sim$ stands for the isomorphic relation and $N(R)$ is the number of graphs isomorphic to $R$. 
Let us denote $\phi_j$ the probability $P$ of pattern $R_j$ given in Figure \ref{Fig:Motifs} as defined in \cite{BCL11}: $\phi_j = P(R_j)$.

\begin{theorem}[\cite{BCL11}] \label{Thm:Bickel} 
Consider a set of fixed patterns ${\bf R} = (R_1,\ldots,R_k)$ with respective sizes $m_j\leq m$ and $\int\int \left(\Phi(u, v)/\phi_1\right)^{2|\Ecal_{R_j}|}\dd u \dd v <\infty$ ($\Phi=\Phi(n)$ and $\phi_1=\phi_1(n)$). Suppose that $(n-1)\phi_1$ is of order $n^{1-2/p}$ or higher. Then,
\[
\sqrt{n}\left((\tilde{P}(R_1),\ldots,\tilde{P}(R_k))-(\Esp\tilde{P}(R_1),\ldots,\Esp\tilde{P}(R_k))\right)\overset{D}{\longrightarrow}  \Ncal(\bf 0, \Sigma_{\bf R}),
\]
where $\tilde{P}(R_j)=\hat{\phi}_1^{-|\Ecal_{R_j}|}\hat{P}(R_j)$ with $\hat{\phi}_1 = \sum_i D_i / [n(n-1)]$. We further have
$\hat{\phi}_1/\phi_1\to^P 1$.
\end{theorem}

\begin{theorem} \label{Thm:W0AsympNormPhi}
Under model $\EGphi$, the statistic $\Wphio$ is asymptotically normal : 
\[(\Wphio - \Esp_\Phi \Wphio)/\Sd_\Phi \Wphio \overset{D}{\longrightarrow}  \Ncal(0, 1),\]
with moments 
\begin{eqnarray}
\Esp_\Phi \Wphio&=&
n^{-1}\left\{n(n-1)^2(\phioo)^2+[1-2(n-1)\phioo]n_1\phi_1
+n_2\phi_2\right\},\label{Thm:W0AsympNormPhi:Esp}\\
  \Var_\Phi \Wphio 
&=&n^{-2} \left\{ 4[1-2(n-1)\phioo]^2 \left(\frac{n_1}2 \phi_1 + n_2\phi_2 + \frac{n_3}4 \phi_1^2 - \frac{n_1^2}4 \phi_1^2  \right) \right. \nonumber\\
&& + 8[1-2(n-1)\phioo] \left[\frac{n_2}2 (2\phi_2 + \phi_3) + \frac{n_3}2 (\phi_5 + 2\phi_6) 
+ \frac{n_4}2 \phi_1\phi_2 - \frac{n_1n_2}4 \phi_1\phi_2 \right] \nonumber\\
&& \left. + 4 \left[\frac{n_2}6 (3\phi_2 + 6\phi_3) + \frac{n_3}2 (4\phi_4 + 2\phi_5 + 2\phi_6 + \phi_7) \right. \right. \nonumber\\
&& \quad \left. \left. + \frac{n_4}4 (4\phi_8 + \phi_9 + 4\phi_{10}) + \left(\frac{n_5}5 - \frac{n_2^2}4\right) \phi_2^2
\right] \right\}\label{Thm:W0AsympNormPhi:Var},
 \end{eqnarray}
where $n_j=\prod_{k=0}^j(n-k)$ and $\phi_j= P(R_j)$ defined just above.
\end{theorem}


\proofbegin 
 The proof relies on the fact that the statistic $\Wphio$ is a linear combination of the $\hat{P}(R_j)$ of three particular patterns $R_j$ to which we will apply Theorem \ref{Thm:Bickel}. Let us begin with the calculation of the moments of $\Wphio$. 
First observe that,
\begin{eqnarray*}
\sum_i [D_i -(n-1)\phioo]^2
&=&n(n-1)^2(\phioo)^2+2[1-2(n-1)\phioo]\sum_{i<j} \Yij\\
&&+2\sum_{1\leq i<j<k\leq n} \Yij\Yik+\Yji\Yjk+\Yki\Ykj\\
&=&n(n-1)^2(\phioo)^2+2[1-2(n-1)\phioo]M_1+2M_2,
\end{eqnarray*}
where 
$$
M_1 = \sum_{1 \leq i < j \leq n} \Yij, \qquad
M_2 = \sum_{1 \leq i < j < k \leq n} \Yij\Yik + \Yij\Yjk + \Yik\Yjk.
$$
Then, we see that,
\begin{eqnarray*}
\Esp_{\Phi} M_1
=\frac{n_1}{2}\phi_1\quad\mbox{and}\quad
\Esp_{\Phi} M_2
=\frac{n_2}{2}\phi_2,
\end{eqnarray*}
which gives $\Esp_\Phi \Wphio$. \\
Next, we calculate the three forthcoming expectations (calculation details are given in Appendix \ref{Proof:Thm:W0AsympNormPhi}):
\begin{eqnarray*}
\Esp_{\Phi} (M_1^2)
&=&\frac{n_1}{2}\phi_1+n_2\phi_2+\frac{1}{4}n_3(\phi_1)^2, \\
\Esp_{\Phi} (M_1M_2)
&=& {\frac{n_2}{2}(2\phi_2+\phi_3) 
+ \frac{n_3}{2}(\phi_5 + 2\phi_6) 
+ \frac{n_4}{4}\phi_1\phi_2} \\
\Esp_{\Phi} (M_2^2)
&=&{\frac{n_2}{6}(3\phi_2 +  6\phi_3) 
+ \frac{n_3}{2}(4\phi_4 + 2\phi_5 + 2\phi_6 + \phi_7) 
+ \frac{n_4}{4}(4\phi_8 + 4\phi_{10} + \phi_9) 
+ \frac{n_5}{4}\phi_2^2 },
\end{eqnarray*} 
which give $\Var_\Phi \Wphio$.

We now turn to the asymptotic normality of $\Wphio $. By using definition \eqref{phat} of $\hat{P}$ and 
the one of $\tilde{P}$ given in Theorem \ref{Thm:Bickel}, we observe that,  
\begin{eqnarray}\label{combi:lin0}
n\Wphio &=& \sum_i [D_i -(n-1)\phioo]^2\nonumber\\
&=&n(n-1)^2(\phioo)^2+[1-2(n-1)\phioo]\sum_{j\neq i} \Yij+\sum_{k\neq j\neq i} \Yij\Yik(1-\Yjk)
+\sum_{k\neq j\neq i} \Yij\Yik\Yjk\nonumber\\
&=&n(n-1)^2(\phioo)^2+[1-2(n-1)\phioo]n_1\hat{P}(R_1)+\frac{1}{3}n_2\hat{P}(R_2)+n_2\hat{P}(R_3)\nonumber\\
&=&n(n-1)^2(\phioo)^2+[1-2(n-1)\phioo]n_1\hat{\phi_1}\tilde{P}(R_1)+\frac{1}{3}n_2\hat{\phi_1}^2\tilde{P}(R_2)
+n_2\hat{\phi_1}^3\tilde{P}(R_3),
\end{eqnarray}
where $R_1$, $R_2$ and $R_3$ are depicted in Figure \ref{Fig:Motifs}. 
Thus we obtain the following linear combination of $\tilde{P}(R_1)$, $\tilde{P}(R_2)$ and $\tilde{P}(R_3)$ :
\begin{eqnarray}\label{combi:lin}
n^{-3/2}\left(\Wphio-(n-1)^2(\phioo)^2\right)
&=&\Theta(n^{1/2})\times\hat{\phi_1}\tilde{P}(R_1)+\Theta(n^{1/2})\times\hat{\phi_1}^2\tilde{P}(R_2)\nonumber\\
&&+\Theta(n^{1/2})\times\hat{\phi_1}^3\tilde{P}(R_3).
\end{eqnarray}
Let us apply the asymptotic normality result of Theorem \ref{Thm:Bickel} to the right-hand side of Equation \eqref{combi:lin}. Since $\hat{\phi}_1\to^P \phi_1$ by Theorem \ref{Thm:Bickel}, we conclude by the Slutsky's lemma. Note that condition $\int\int \left(\Phi(u, v)/\phi_1\right)^{2|\Ecal_{R_j}|}\dd u \dd v <\infty, j=2,3$ is fulfilled because $\Phi\leq 1$ and $\phi_1$ are constants.
\proofend

\begin{remark} \label{Rem:LargerVariance}
The test statistics $\Wo$ in the independent case and $\Wphio$ in the exchangeable case measure both the discrepancy between the
observed degrees and their expected values under specifics models. 
Let us stress that the latent layer in the exchangeable case implies an additional variability of the degrees. The third term in Equation \eqref{combi:lin} is a consequence of this additional variability.
\end{remark}

\subsubsection*{Particular cases: SBM and EDD}\label{Sec:SBM:EDD}
Because SBM and EDD are special cases of the $W$-graph, all results above apply to them. Interestingly, for both models, the critical calculation of coefficients $\phi_1$ to $\phi_{10}$ can be achieved exactly. Indeed, the calculation of the first two moments of pattern counts under SBM and EDD is explicitly addressed in \cite{PDK08}. In this reference, it is already observed that patterns 4 to 10 from Figure \ref{Fig:Motifs} need to be considered as 'super-patterns' (or 'super-motifs') of patterns 2 and 3 and that the variance of the count of a given pattern depends on the expected frequency of its super-patterns. \\
The formula of $\phi_j$ for SBM is explicitly in \cite{PDK08}. Denoting $\alpha_k$ the probability for any given node to belong to group $k$ ($1 \leq k \leq K$), we have that
$$
\phi_j = P(R_j) = \sum_{k_1}^K \dots \sum_{k_{p_j}}^K \alpha_{k_1} \dots \alpha_{k_{p_j}} \prod_{1
\leq u < v \leq p_j} \pi_{k_u k_v}^{m^j_{uv}}
$$
where $p_j$ stands for number of nodes in pattern $R_j$ and $m^j_{uv}$ is 1 if nodes $u$ and $v$ are
connected in pattern $R_j$ and 0 otherwise. \\
The EDD model is also studied in \cite{PDK08} but needs to be adapted to the $W$-graph framework. For $\phi(u, v) = g(u) g(v)$, we have that
$$
\phi_j = \prod_{u=1}^{p_j} g_{d^j_u}, 
\qquad \text{where} \quad
g_k = \int_0^1 g^k(u) \dd u
$$
and $d^j_u$ stands for the degree of node $u$ within the pattern $R_j$. Some examples are
$$
\phi_1 = g^2_1, \qquad
\phi_2 = g^2_1 g_2, \qquad
\phi_3 = g^3_2, \qquad
\phi_4 = g_1 g^2_2 g_3, \qquad
\phi_{10} = g^3_1 g_2 g_3.
$$
\paragraph{Plug-in version of the test.}
In many situations, $\phio$ is unknown and one needs to resort to an estimate $\phioh$. The question is then to understand the asymptotic behaviour of $(\Wphiho - \Esp_\Phi \Wphiho)/\Sd_\Phi \Wphiho$. A first strategy would consist in estimating $\phio$ from the data. Still, few results are available regarding the statistical properties of the graphon estimates. More recently, \cite{GaL17b} considered a simpler statistic, the moments of which can be estimated via the empirical counts of the  patterns $R_1, R_2, R_3$ (see Figure \ref{Fig:Motifs}). They proved the asymptotic normality of its plug-in version in the degree-corrected SBM model. In our case, this would require to establish asymptotic results about quantities that combine patterns $R_1$ to $R_{10}$, in a particularly intricate manner.

\subsection{Test and power} \label{Sec:TestDMS:Graphon}

We now study the test of $H_0=\EGphio$ versus $H_1 = \EGphi$. The next Corollaries provide the null distribution of the test statistic $\Wphio$ and the power of the associated test. They are direct consequences of Theorem \ref{Thm:W0AsympNormPhi}.

\begin{corollary} \label{Cor:W0AsympNormPhi0}
Under the model based on $\phio$ the statistic $\Wphio$ is asymptotically normal with moments expressed as those of Theorem \ref{Thm:W0AsympNormPhi} with all $\phi_j$ replaced by $\phi_j^0$.
\end{corollary}
Recall that the particular terms $\dij = \pij - \poij$ appear in the moments of $\Wo$ under model $HER(\pg)$ whereas it is not the case anymore under $\HERpgo$ (see Theorem \ref{Thm:W0AsympNormHER} and Corollary \ref{Cor:W0AsympNormHER0} in sections \ref{Sec:MeanSqDeg} and \ref{Sec:TestDMS}). Notice that this simple measure of discrepancy between two alternative models is not visible in the moments of $\Wphio$ but spread out all differences between $\phi_j$ and $\phi^0_j$. 

A formal test with asymptotic level $\alpha$ can be constructed based on Corollary \ref{Cor:W0AsympNormPhi0}, which reject $H_0$ as soon as $\Wphio$ exceeds $\Esp_\phio \Wphio + t_{1-\alpha} \Sd_\phio \Wphio$. The expression of its power follows.

\begin{corollary}
The asymptotic power of the considered test is 
\begin{align}
 \pi(\pg) = 1 - \Phi\left(\left(\Esp_\phio \Wphio + t_{1-\alpha}\Sd_\phio \Wphio - \Esp_\Phi \Wphio\right) \left/ \Sd_\Phi \Wphio \right.\right).\label{Cor:EG:Power}
\end{align}
\end{corollary}

\begin{remark} 
Let consider the test of $H_0 = ER$ versus $H_1 = \EGphi$. This simply corresponds to the degree variance test based on the statistic $V$ described in Section \ref{Sec:TestDV}.
\end{remark} 

The following corollary gives a sufficient condition on the departure between $\Phi$ and $\Phi^0$ to ensure that the proposed test is asymptotically powerful.
\begin{corollary}\label{Cor:W0Power}
For functions $\Phi^0$ and $\Phi$, define
$$
\Delta_n(\Phi^0, \Phi) :=\Esp_{\Phi} \Wphio - \Esp_{\phio} \Wphio= [1-2(n-1)\phioo]n_1(\phi_1-\phi_1^0) +n_2(\phi_2-\phi_2^0).
$$
If $\Delta_n(\Phi^0, \Phi) > 0$, then the test $H_0 = \EGphio$ versus $H_1 = \EGphi$ is asymptotically powerful.
\end{corollary}
\proofbegin
The proof follows the line of Corollary \ref{Cor:W0AsympNormHER0:Pwer}. The expression of $\Delta_n(\Phi^0, \Phi)$ comes from \eqref{Thm:W0AsympNormPhi:Esp}. Because functions $\Phi$ and $\Phi^0$ are fixed, we have that $\Delta_n(\Phi^0, \Phi) = \Theta(n^3)$. Furthermore, from \eqref{Thm:W0AsympNormPhi:Var}, we have that $\Sd_{\Phi} \Wphio = \Theta(n^{3/2})$ and $\Sd_{\phio} W_\phio = \Theta(n^{3/2})$. As a consequence, if $\Delta_n(\Phi^0, \Phi) > 0$, the negative argument of $\Phi$ in \eqref{Cor:HER:Power} goes to infinity at rate $n^{3/2}$, which concludes the proof.
\proofend

\begin{remark} In Corollary \ref{Cor:W0Power}, \\
$(i)$ the condition only depends on the relative frequencies of $R_1$ and $R_2$. If $\Phi$ and $\phio$ have the same $\phi_1$ and $\phi_2$ but differ in terms of, say $\phi_k$ ($k > 2$) the proposed test may no be able to detect the discrepancy.\\
$(ii)$ observe that, if the function $\Phi$ depends on $n$ (and is denoted $\Phi_n$), the asymptotic power is still guaranteed as long as $\Delta_n(\Phi^0, \Phi_n) = \Theta(n^\alpha) > 0$ with $\alpha > 3/2$.
\end{remark}

\subsubsection*{Illustration}
As  an illustration  of the  proposed test,  we consider  the networks
described in  Section \ref{Sec:Illustration}. The question  is to know
if a  fitted graphon is sufficient  to explain the heterogeneity  of a
network, at least  in terms of degrees. To address  this question, for
each network separately, we estimated a graphon function using the variational
  expectation maximization of \cite{DPR08} to provide estimates of the
  SBM model parameters and build the corresponding block-wise constant
  graphon function. The number of blocks was estimated using the model
  selection criterion considered in \cite{DPR08}. This is implemented in the package {\tt mixer} (available on the \url{https://cran.r-project.org/}). We then calculated the moments of the graphon and applied the degree mean square test to check if the fitted graphon is sufficient to explain the heterogeneity of the network. The results are given in Table \ref{Tab3:EG-SBM}.

\begin{table}[!ht]
  \begin{center}
      \caption{Degree mean square EG test for an SBM-graphon. TestStat $= ({\Wphiho-\Esp_{\hat{\Phi}^0}})/{\Sd_{\hat{\Phi}^0}}$. \label{Tab3:EG-SBM}}
      
    \begin{tabular}{lccccccc}
 Network & $n$ & density & $K$ & $\Wphiho$ & $\Esp_{\hat{\Phi}^0}$ & $\Sd_{\hat{\Phi}^0}$ & TestStat \\ \hline
 Karate & 34 & 0.139 & 4 & 14.6 & 15.57 & 6.16 & -0.16 \\ 
Tree & 51 & 0.54 & 5 & 163.14 & 162.84 & 17.31 & 0.02 \\ 
Fungi & 154 & 0.227 & 15 & 597.6 & 584.42 & 116.63 & 0.11 \\ 
Blog & 196 & 0.075 & 11 & 104.72 & 92.77 & 25.89 & 0.46 \\ 
CKM & 219 & 0.015 & 3 & 3.9 & 4.04 & 0.76 & -0.18 \\ 
FauxDixon & 248 & 0.02 & 5 & 16.78 & 11.97 & 1.94 & 2.48 \\ 
AdHealth & 530 & 0.007 & 4 & 10.7 & 7.54 & 1.42 & 2.22 \\ 
    \end{tabular}
  \end{center}
\end{table}

{Using the normal approximation for the distribution of $\Wphiho$ under $H_0$, the $\EGphioh$ model is rejected for two of these networks: FauxDixon and AdHealth. The highest test statistic is observed for the FauxDixon network, which has actually been simulated under a model that does not belong to the class of $\EGphi$.}


\subsection{Case of sparse graphs} \label{Sec:Sparse:Graphon}

\noindent The following theorem discusses the validity of Theorem \ref{Thm:W0AsympNormPhi} when considering sparse graphs, namely when $\phi_1=\phi_1(n)$ vanishes as $n$ grows with a rate we specify.
\begin{proposition}
Under the model based on the graphon $\Phi$ such that $\phi_1$ and $\phioo$ are of order $n^{-2/3}$ or higher, if $\iint \left(\Phi(u, v)/\phi_1(n)\right)^{6}\dd u \dd v <\infty$ then the statistic $\Wphio$ is asymptotically normal. 
\end{proposition}

\proofbegin 
We apply Theorem \ref{Thm:Bickel} to a function of $\Wphio$ which is a linear combination of $\tilde{P}(R_1)$, $\tilde{P}(R_2)$ and $\tilde{P}(R_3)$ involving the quantity $\hat{\phi_1}/\phi_1$, where $R_1$, $R_2$ and $R_3$ refer to the patterns from Figure \ref{Fig:Motifs}. Equations \eqref{combi:lin0}--\eqref{combi:lin} state that
\begin{eqnarray*}
n^{-3/2}\left(\Wphio-(n-1)^2(\phioo)^2\right)&=&\Theta(\phioo\phi_1)\times\Theta(n^{1/2})\times\frac{\hat{\phi_1}}{\phi_1}\tilde{P}(R_1)\nonumber\\
&&+\Theta(\phi_1^2)\times\Theta(n^{1/2})\times \left(\frac{\hat{\phi_1}}{\phi_1}\right)^2\tilde{P}(R_2)\nonumber\\
&&+\Theta(\phi_1^3)\times \Theta(n^{1/2})\times \left(\frac{\hat{\phi_1}}{\phi_1}\right)^3\tilde{P}(R_3).
\end{eqnarray*}
The asymptotic normality of $\sqrt{n}\left(\tilde{P}(R_1),\tilde{P}(R_2),\tilde{P}(R_3)\right)$ holds under conditions : \\
$\int\int \left(\Phi(u, v)/\phi_1\right)^{2|\Ecal_{R_j}|}\dd u \dd v <\infty$ with $|\Ecal_{R_j}|\leq3$ and $\phi_1$ being of order $n^{-2/p}$ or higher with $p=3$. Now, we observe that under the condition that $\phi_1$ and $\phioo$ are of order $n^{-\alpha}$ for $0<\alpha<2/3$,
\begin{eqnarray*}\label{combi:lin:sparse}
n^{-3/2+2\alpha}\left(\Wphio-(n-1)^2(\phioo)^2\right)&=&\Theta(n^{1/2})\times \frac{\hat{\phi_1}}{\phi_1}\tilde{P}(R_1)
+\Theta(n^{1/2})\times \left(\frac{\hat{\phi_1}}{\phi_1}\right)^2\tilde{P}(R_2)\nonumber\\
&&+\Theta(n^{1/2-\alpha})\times \left(\frac{\hat{\phi_1}}{\phi_1}\right)^3\tilde{P}(R_3).
\end{eqnarray*}
Since $\hat{\phi}_1/ \phi_1\to^P1$ by Theorem \ref{Thm:Bickel}, 
we conclude by applying the asymptotic normality result of the same theorem to the right-hand side of the equality above combined with the Slutsky's lemma. Note that the third term mentioned in Remark \ref{Rem:LargerVariance} is negligible.
\proofend

\section{Simulation study} \label{Sec:Simus}

We designed a simulation study to  assess the performance of the tests described above. More specifically, our purpose is to evaluate the power of these tests for various graph sizes and  densities (mean connectivities). We also aim at illustrating for which graph size the asymptotic normal approximation is accurate; we especially focus on this point in the sparse regime. 

\subsection{Simulation design}

\paragraph{Design for the independent case.}
We designed our simulation so that to mimic the situation where an heterogeneous model $\HERpgo$ is considered, which still misses some heterogeneity. More specifically, each node $i$ was associated with a vector of covariates $x_i \in \mathbb{R}^d$ (all values were drawn i.i.d. with standard Gaussian distribution and $d$ was set to 3). Each edge $(i, j)$ was then associated with the covariate vector $\xij = \sqrt{\pi}|x_i - x_j|/2$ so that all $x_{ij}$  are positive with mean 1. The edges were then drawn according a logistic model: $\text{logit}(\pij) = a + \xpij \beta_1$ where $\beta_1^\intercal = {[\beta_0^\intercal \; \beta]}^\intercal \in \mathbb{R}^d$, $\beta_0 \in \mathbb{R}^{d-1}$. The constant $a$ was set to preserve the mean connectivity, denoted $\rho^{*}$ in the sequel. 
The probability matrix $\pgo = [\poij]$ of the null model was defined according to the same logistic model, removing the last covariate, namely $\text{logit}(\poij) = a_0 + \xopij\beta_0$, where $\xoij$ is $\xij$ deprived from its last coordinate. Hence, the discrepancy between the null hypothesis and the true model is measured by the coefficient $\beta$ of the last covariate. All $\beta_0$'s were set to $1$ except $\beta$ which ranged from 0 to 2. We also studied the behaviour of the plug-in version  as defined in the paragraph 'Plug-in version of the test' at the end of Section \ref{Sec:MeanSqDeg}. 

\paragraph{Design for the exchangeable case.} 
We designed a situation where a null block-wise constant graphon $\Phi^0$, associated to a SBM model, is contaminated by an alternative graphon of the form considered in \cite{LaR15}. Thus, graphs were sampled from an $\EGphi$ model where $\Phi(u, v) = \Phi^0(u, v)\rho \beta^{2}u^{\beta-1}v^{\beta-1}$. Note that $\Phi$ induces a random graph model related to the degree corrected SBM model of \cite{karrer11} which has received strong attention in the last five years. This model, by characterizing explicitly the degrees of the vertices, is often employed as an alternative to the standard SBM model. Note however that in its original form the degree corrected SBM model is not exchangeable since the degree parameters are fixed. Conversely, $\Phi$ induces an exchangeable model here since the degree terms $u^{\beta-1}$ and $v^{\beta-1}$ are random.   For the null graphon $\Phi^0$, we considered a SBM with 2 blocks, with the same proportions. Moreover, $\Phi^0$  was given a product form such that  $\Phi^0(u, v)=\eta_k\eta_\ell$ if $u$ and $v$ are  in blocks $k$ and $\ell$, respectively. We set $\eta_1=0.4$ and $\eta_2=0.5$. In this simulation framework, the discrepancy between the null hypothesis and the true model is measured by the term $\beta$ which ranges from 1 to 2 and controls the imbalance of the expected degrees of the nodes. The null graphon is retrieved when $\beta = 1$. Finally, the term $\rho$ was set in order to obtain the desired mean connectivity $\rho^{*}$. 

Note that is both designs, the density of the network is kept constant equal to {$\rho^{*}$} when going away from the null model. Therefore, the departure from $H_0$ detected by the tests is not due to a mean degree difference. In both designs, $\beta$ measures the departure from the null model, although the its nominal values are not comparable from one design to another. 
1\,000 simulations were ran for each combination of the parameters $(n, \rho^{*}, \beta)$.

\paragraph{Sparse graphs.}
For both tests, we considered sparse graphs in the setting described in Sections \ref{Sec:Sparse} and \ref{Sec:Sparse:Graphon}. We focused on the asymptotic normality of the degree mean square statistic under the null hypothesis. To this aim, we designed  a reference null probability matrix $\pgso$ and a reference null graphon $\phiso$ as described above. We then considered the two sparsity scenarios:
\begin{itemize}
 \item vanishing connection probabilities: $\poij = \psoij n^{-a}$ and $\Phi(u, v) = n^{-a} \phiso(u, v)$;
 \item sparse connection probabilities: $\poij = \psoij$ with probability $n^{-b}$ and 0 otherwise.
\end{itemize}
The second scenario does not make sense for the $\EGphi$ test. {The mean connectivity $\rho^{*}$ was set to 0.1}. The density of the graphs therefore decrease as $\rho^{*} n^{-a}$ and $\rho^{*} n^{-b}$, respectively.

\paragraph{Criteria.}
For each parameter configuration, we computed the moments of the respective statistics and derived the theoretical power. Based on the replicates, we estimated the empirical power, and its plug-in version in the independent case. For the sparse setting, the proximity with the normal distribution was investigated plotting the empirical quantiles versus the theoretical Gaussian quantiles (QQ-plots).

\subsection{Results}

\begin{figure}[!ht]
 \begin{center}
 \includegraphics[scale=0.85]{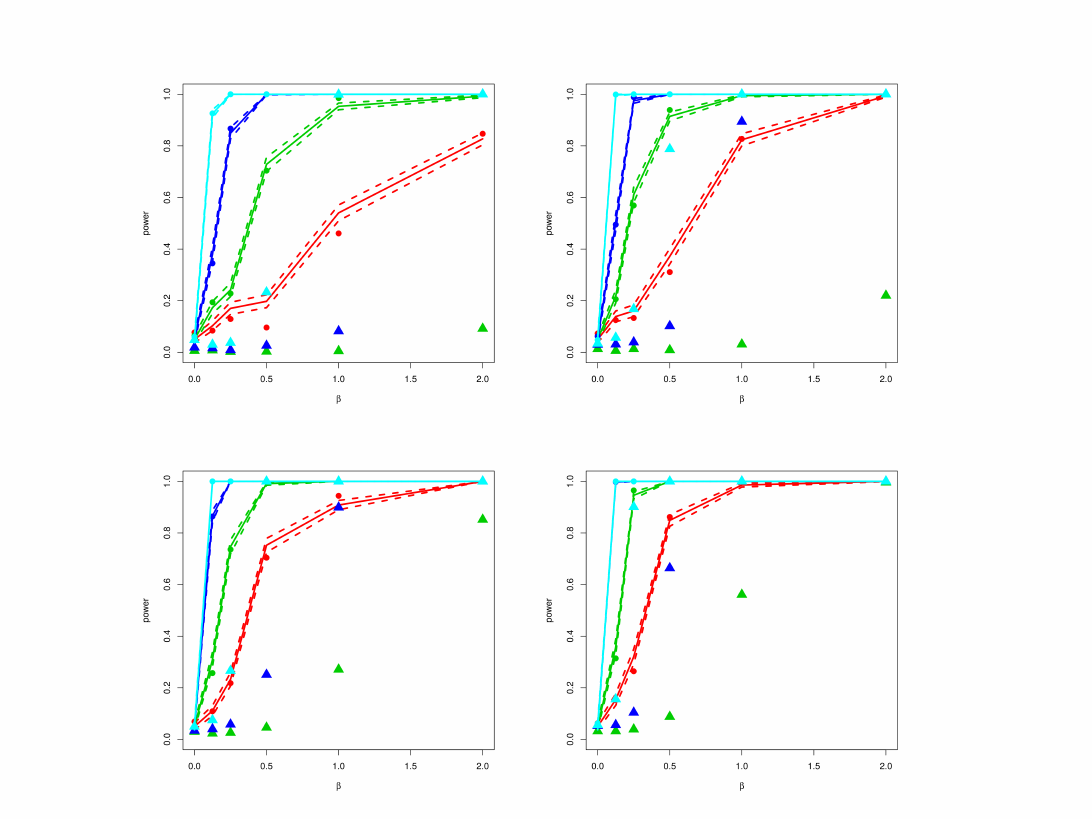}
  \caption{Power of the degree mean square test in the HER design, as a function of {$\beta$} the
    effect of the last covariate. 
    Top left ($\log_{10} \rho^{*}= -2.5$), top right ($\log_{10} \rho^{*}= -2$), bottom left ($\log_{10} \rho^{*}= -1.5$), bottom right ($\log_{10} \rho^{*}= -1$). 
    Color refers to the graph size: $n = 32$ (red), $100$ (green), $316$ (blue), $1\,000$ (cyan) (green, blue and cyan curves and points overlap in the last panels). 
    Points = empirical power (average on $1\,000$ simulations): 
    dots = $\Wo$ test,
    solid line = theoretical power,
    dashed line = binomial confidence interval for $1\,000$ simulations, 
    triangles = $\Who$ test (for $n \geq 100$).
    \label{Fig:PowerW} }
 \end{center} 
\end{figure}

\begin{figure}[!ht]
 \begin{center}
 \includegraphics[scale=0.85]{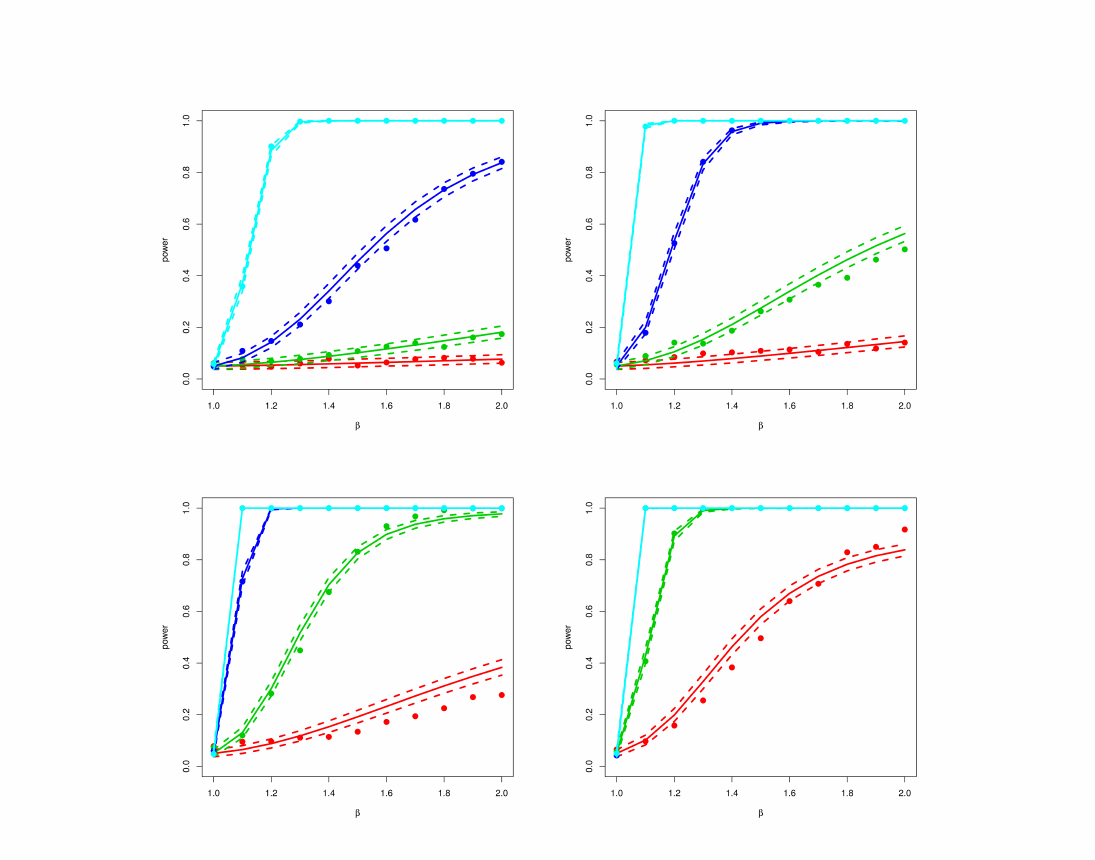}
  \caption{Power of the degree mean square test in the EG design, as a function of $\beta$ which controls the degree imbalance. Same legend as in Figure \ref{Fig:PowerW}.
  \label{Fig:PowerW_exchangeable}}
 \end{center}
\end{figure}

\paragraph{Power and asymptotic normality.}
The power curves of the degree mean square tests in the independent and exchangeable cases are given in 
s \ref{Fig:PowerW} and \ref{Fig:PowerW_exchangeable}, respectively. As expected, the power increases with the departure $\beta$, the graph size $n$ and the network density $\rho^{*}$. We remind that the departure parameter $\beta$ can not be compared between the two figures.
The binomial confidence interval around the theoretical power informs us about the convergence to the asymptotic normality. We observe that the empirical power (dots) falls within this interval showing that the normal approximation is accurate for reasonably large ($n > 100$) graphs. This does not hold for the empirical version of the HER test (triangles), which suggests that the cumulative effect of all the estimation errors $| \phoij - \poij |$ on $\Who$ vanishes later than the convergence of $\Wo$ to normality.
The power of both tests also depends on the density of the graph; it is satisfying for $\rho^* \geq 1\%$ in the independent case and for $\rho^* \geq 3\%$, in the exchangeable case. As for the empirical version of the HER test,  it becomes reasonable only when $n$ reaches $300$, whatever the density.

\begin{figure}[!h]
\begin{center}
\includegraphics[scale=0.85]{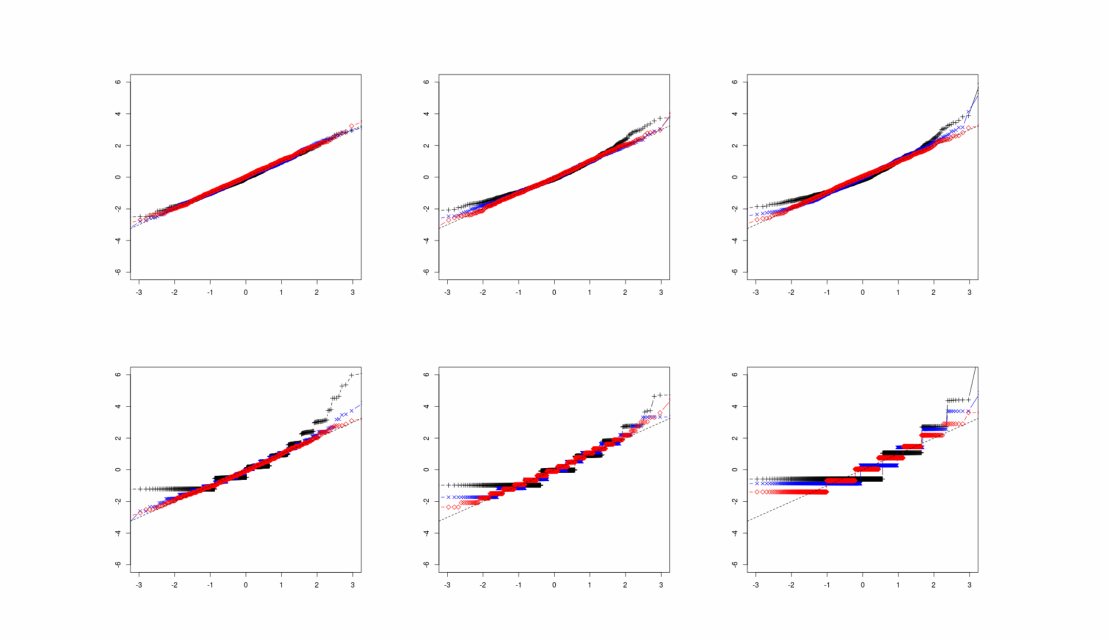}
  \caption{QQ-plots of the degree mean square statistics $\Wo$ in the HER design, for vanishing connection probabilities: $\pij = \psij n^{-a}$ and initial mean density $\rho^{*} = 0.1$.
  From top left to bottom right: $a = 0, 0.4, 0.8, 1.2, 1.4, 1.6$.
  Graph size $n$ = 100 ($+$), 1\,000 ($\textcolor{blue}{\times}$) and 10\,000 ($\textcolor{red}{\diamond}$). 
  \label{Fig:AsNormWsmallProb}}
 \end{center}
\end{figure}

\begin{figure}[!h]
 \begin{center}
 \includegraphics[scale=0.85]{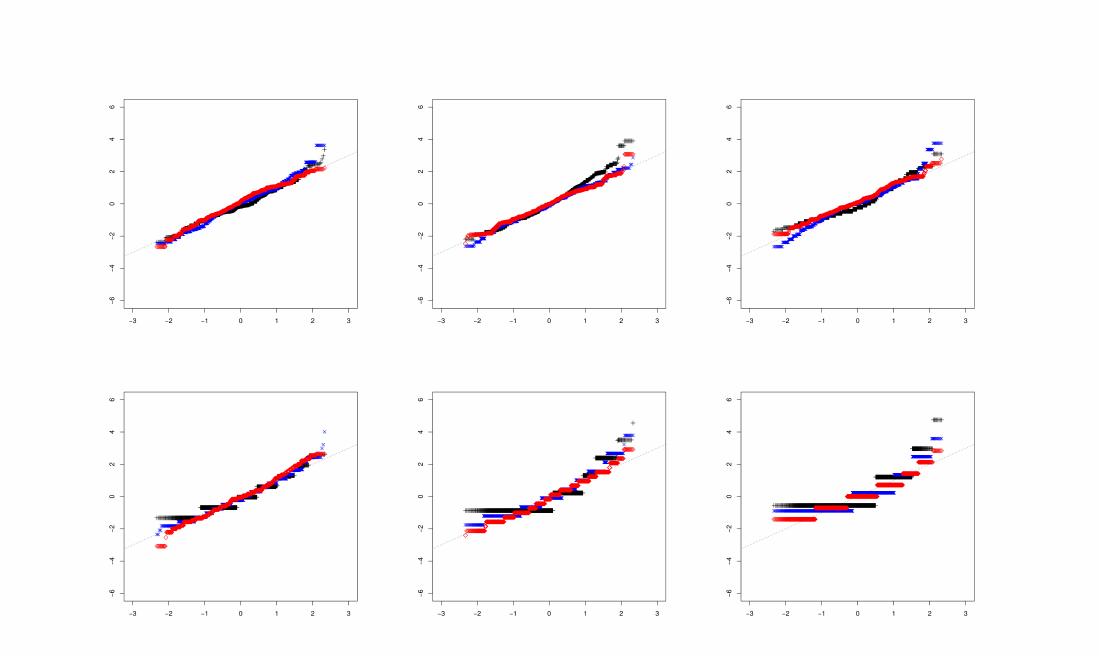}
  \caption{QQ-plots of the degree mean square statistics $W_{\Phi^{0}}$ in the EG design, for vanishing connection probabilities: $\Phi(u, v) = n^{-a}\Phi^{0*}(u, v)$ and initial mean density $\rho^{*}=0.1$. 
  Same legend as in Figure \ref{Fig:AsNormWsmallProb}.
  \label{Fig:AsNormWexchangeableSparse}}
 \end{center}
\end{figure}

\paragraph{Sparse graphs.} 
Figures \ref{Fig:AsNormWsmallProb} and \ref{Fig:AsNormWexchangeableSparse} display the QQ-plots of the standardized $\Wo$ and $\Wphio$ statistics under the vanishing probabilities scenario for graphs with several sizes. Remember that the larger the power $a$, the sparser the graph.
We observe again that normality holds for the non sparse graphs ($a = 0$) even for $n = 100$, but the departure is visible for $n = 100$ as soon as $a \geq 0.4$. The same is observed for $n = 1\,000$, although a bit later ($a \geq 0.8$). For the largest graph ($n = 10\,000$), normality holds until $a \simeq 1.2-1.4$ but does not seem to be reached for higher sparsity regimes. 
As expected, in the very sparse regime, normality can only be relied on for very large graphs.
Similar conclusions can be drawn for the sparse probabilities scenario, each distribution being slightly closer to normal. 

\section*{Acknowledgements} 
This work has been partially funded by the research grant NGB (ANR-17-CE32-0011).
We thank Pr Bloznelis for providing us with his report referred to as \cite{Blo05}.

\bibliography{biblio}
\bibliographystyle{chicago}

\appendix
\section{Appendix}


\subsection{Proof of Corollary \ref{Thm:VhAsympNormHER}} \label{Proof:Thm:VhAsympNormHER}

Let express $\Vh$ as follows.
\begin{eqnarray}\label{wn}
n^2\Vh&=&\frac1{2} \sum_{i \neq j} \left(D_i - D_j\right)^2\nonumber\\
&=&2 (n-2) \sum_{1\leq i<j\leq n} \Yij\nonumber\\
&& + 2(n-4)\sum_{1\leq i<j<k\leq n} \left\{\Yij\Yik+\Yij\Yjk+\Yik\Yjk\right\}\nonumber\\
&& - 8\sum_{1\leq i<j<k<l\leq n} \left\{\Yij\Ykl+\Yik\Yjl+\Yil\Yjk\right\}.
\end{eqnarray}
Then we write the Hoeffding decomposition of $\Vh$  :
\begin{eqnarray}\label{hoef}
\Vh&=&P_{\emptyset}\Vh+\sum_{1\leq i<j\leq n} P_{\{ij\}}\Vh+
\sum_{1\leq i<j<k\leq n}\left\{P_{\{ij,ik\}}\Vh+P_{\{ij,jk\}}\Vh+P_{\{ik,kj\}}\Vh\right\}\nonumber\\
&&+\sum_{1\leq i<j<k<l\leq n}\left\{P_{\{ij,kl\}}\Vh+P_{\{ik,jl\}}\Vh+P_{\{il,jk\}}\Vh\right\}.
\end{eqnarray}
Taking all projections with respect to $HER(\pg)$, we have
\begin{eqnarray}
P_{\emptyset}\Vh&=&\frac{1}{2n^2}\left(4(n-2)\sum_{1\leq i<j\leq n} \pij+4(n-4)\sum_{1\leq i<j<k\leq n}\left\{\pij\pik+\pij\pjk+\pik\pjk\right\}\right)\nonumber\\
&&-\frac{8}{n^2}\sum_{1\leq i<j<k<l\leq n}\left\{\pij\pkl+\pik\pjl+\pil\pjk\right\},\nonumber
\end{eqnarray}
which gives the expectation. The other projections provide the variance. We have,
\begin{align}
&P_{\{ij\}}\Vh=\frac{1}{2n^2}\Ytij\left(4(n-2)+4(n-4)\sum_{k\notin(i,j)}(p_{ik}+p_{jk})-16\sum_{k<l\notin(i,j)}\pkl\right),\label{pij}\\
&P_{\{ij,ik\}}\Vh=\frac{2(n-4)}{n^2}\Ytij\Ytik,\qquad\mbox{ and }\qquad P_{\{ij,kl\}}\Vh=-\frac{8}{n^2}\Ytij\Ytkl.\label{pijkl}
\end{align}
So,
\begin{align} 
&n^4 \Var P_{\{ij\}}\Vh=\sigmaij\left(2(n-2)+2(n-4)\sum_{k\notin(i,j)}(p_{i,k}+p_{j,k})-8\sum_{k<l\notin(i,j)}\pkl\right)^2,\label{vpij}\\
&n^4 \Var P_{\{ij,ik\}}\Vh= 4(n-4)^2\sigmaij\sigma^2_{ik},\qquad\mbox{ and }\qquad
n^4 \Var P_{\{ij,kl\}}\Vh=64\sigmaij\sigmakl,\label{vpijkl}
\end{align}
and the variance of $\Vh$ follows by summing over all indexes.\\
As for the asymptotic normality, we consider $\Vh-\Esp\Vh=\Vh^*-\Esp_{\bf{p}}\Vh+\Vh-\Vh^*$, with $\Vh^*=P_{\emptyset}\Vh+\sum_{1\leq i<j\leq n}P_{\{ij\}}\Vh$. In order to show that  that $\Vh^*-\Esp_{\bf{p}}\Vh$ is asymptotically normal, we apply Theorem \ref{Thm:LindebergLevy} to the projections $P_{\{ij\}}\Wo$ (which stand for the $\Xnu$) by using Remark \ref{rem1} and Assumption \ref{Assumption}.
The $a_{n\{ij\}}=\Theta(1)$ expressed in \eqref{pij} stand for $a_{nu}$. Since $B^2_n=\Var_{\bf{p}}\left(\Vh^*-\Esp_{\bf{p}}\Vh\right)=\Theta(n^2)$, we conclude that the Lindeberg condition is fulfilled because, for any $\epsilon$, each $\anu$ becomes smaller than $\epsilon B_n$ when $n$ goes to infinity. Now we consider $\Vh-\Vh^*$ as the linear combination of the projections $P_{\{ij,ik\}}\Vh$ and $P_{\{ij,kl\}}\Vh$. We notice that $a_{n\{ij,ik\}}$ and $a_{n\{ij,kl\}}$ given in \eqref{pijkl} equal $\Theta(n^{-1})$ and $\Theta(n^{-2})$ respectively, and thus that $\Var_{\bf{p}}\left(\Vh-\Vh^*\right)=\Theta(n)$. We conclude to the asymptotic normality of $\Vh$ by combining the one of $\Vh^*-\Esp\Vh$ and the fact that $\Var_{\bf{p}}\left(\Vh-\Vh^*\right)/\Var_{\bf{p}}\Vh^*\to 0$ as $n\to\infty$.

\subsection{Degree variance test power}\label{Lemma:PowerER}
\begin{lemma}
 \label{Cor:VhTestER}
Under model $ER$ and Assumption \ref{Assumption}, the degree variance is asymptotically normal:
 $$
 \left(\Vh - \Esp_{\hat{p}} \Vh\right)  / \Sd_{\hat{p}}\Vh \overset{D}{\longrightarrow} \Ncal(0, 1).
 $$
\end{lemma}
\proofbegin
The proof relies on the concentration of $\ph$ around $p$ and on Slutsky's lemma (see, e.g., Theorem 4.4, p.27 in Billingsley \cite{Bil68}). 
First, write the statistic based on $\Vh$ as
\[\frac{\Vh - \Esp_\ph \Vh}{\Sd_\ph\Vh}
=
\frac{\Sd_p \Vh}{\Sd_\ph \Vh} \left(\frac{\Vh - \Esp_p \Vh}{\Sd_p\Vh} + \frac{\Esp_p \Vh-\Esp_\ph \Vh}{\Sd_p\Vh} \right).
\]
Then note that, under $ER(p)$, $(\ph - p) = \Theta_P(n^{-1})$, so $(\ph\qh - pq) =\Theta_{\mathbb{P}}(n^{-2})$, where $\qh$ stands for $1- \ph$. 
According to the moments given in Corollary \ref{Thm:VhAsympNormHER}, we have that $\Esp_p \Vh =\Theta(n)pq$ and $\Var_p \Vh = \Theta(1) pq + \Theta(n) p^2q^2$. This entails that $\Esp_p \Vh-\Esp_\ph \Vh =\Theta_{\mathbb{P}}(n^{-1})$ and 
$\Var_\ph \Vh - \Var_p \Vh =\Theta_{\mathbb{P}}(n^{-2})$, so ${{\Sd_p \Vh}/{\Sd_\ph \Vh}}$ converges in probability to 1 and ${(\Esp_p \Vh-\Esp_\ph \Vh)}/{\Sd_p \Vh}$ converges in probability to 0. The result then follows from Slutsky's lemma, used twice. \proofend

\begin{lemma}
\label{Lem:VphEquivVpb}
We have
 \[\left(\Vh - \Esp_{\hat{p}} \Vh\right)  / \Sd_{\hat{p}}\Vh-
\left(\Vh - \Esp_{\bar{p}} \Vh\right)  / \Sd_{\bar{p}} \Vh  \overset{\mathbb{P}}{\longrightarrow} 0,\] 
where $\pb=[n(n-1)]^{-1}\sumij\pij$.
\end{lemma}
The proof of this Lemma is similar to this of Lemma \ref{Cor:VhTestER} and results from the concentration of $\ph$ around $\pb$.

\subsection{Proof of Corollary \ref{Prop:VhAsympNormHER} \label{Proof:Prop:VhAsympNormHER}}

The proof follows the line of this of Proposition \ref{Prop:W0Sparse} under Assumption \ref{Assumption}. We begin with the asymptotic normality of $\Vh^*-\Esp_{\bf{p}}\Vh$. Since $\sum_{k\notin(i,j)} \pik=\Theta(n^{1-a-b})$ and $\sum_{k<l\notin(i,j)}\pkl=\Theta(n^{2-a-b})$, we see that $a_{n\{ij\}}=\Theta(n^{-(a+b)})$ if $a+b<1$ and $\Theta(n^{-1})$ if $a+b>1$ ($a_{n\{ij\}}$ are given in assertion \eqref{pij}). Therefore, we have $\Var_{\bf{p}} P_{\{ij\}} V= \Theta\left(n^{-3a-2b}\right)$ if $a+b<1$ and $\Theta\left(n^{-a-2}\right)$ if $a+b>1$. Combining this with the number of non-zero terms which is $\Theta(n^{2-b})$, we get that $B_n^2 =\Theta\left(n^{2-3(a+b)}\right)$ if $a+b<1$ and $\Theta\left(n^{-(a+b)}\right)$ if $a+b>1$.
Comparing $A_n^2(\epsilon)$ with $B_n^2$, we see that the Lindeberg condition is fulfilled for $a+b < 2$.\\
Now we consider $\Vh-\Vh^*$ as the linear combination of the projections $P_{\{ij,ik\}}\Vh$ and $P_{\{ij,kl\}}\Vh$. We see that $a_{n\{ij,ik\}}=\Theta(n^{-1})$ and $a_{n\{ij,kl\}}=\Theta(n^{-2})$ ($a_{n\{ij,ik\}}$ and $a_{n\{ij,kl\}}$ are given in assertion \eqref{pijkl}). Therefore, we have 
$\Var_{\bf{p}} P_{\{ij,ik\}} \Vh =\Theta\left(n^{-2a-2}\right)$ and $\Var_{\bf{p}} P_{\{ij,kl\}} \Vh = \Theta\left(n^{-2a-4}\right)$. Since the number of non-zero terms in the sums is $\Theta(n^{3-2b})$ and $\Theta(n^{4-2b})$ respectively, we have therefore $\Var_{\bf{p}}\left(\Wo-\Wo^*\right)=\Theta(n^{-2(a+b)+1})$.\\
We conclude to the asymptotic normality of $\Vh$ by combining the one of $\Vh^*-\Esp_{\bf{p}}\Vh$ under condition $a+b < 2$ and the fact that $\Var_{\bf{p}}\left(\Vh-\Vh^*\right)/\Var_{\bf{p}}\Vh^*\to 0$ as $n\to\infty$.

\subsection{Moments of $\Wphio$ in the proof of Theorem \ref{Thm:W0AsympNormPhi} \label{Proof:Thm:W0AsympNormPhi}}

We have
\begin{eqnarray*}
\Esp_{\Phi}(M_1^2)
&=&\Esp_{\Phi}\left(\sum_{i<j}\Yij^2\right)
+2\Esp_{\Phi}\left(\sum_{1\leq i<j<k\leq n}\Yij\Yik+\Yji\Yjk+\Yki\Ykj\right)\\
&&+2\Esp{\Phi}\left(\sum_{1\leq i<j<k<l\leq n} \Yij\Ykl+\Yik\Yjl+\Yil\Yjk\right) \\
&=&\frac{n_1}{2}\phi_1+n_2\phi_2+\frac{1}{4}n_3(\phi_1)^2,
\end{eqnarray*}
and
\begin{eqnarray*}
\Esp_{\Phi}(M_1M_2)
%
&=& {\binom {3} {2}}\Esp_{\Phi}\left(\sum_{1\leq i<j<k\leq n}\Yij^2\Yik+\Yij^2\Yjk+\Yij\Yki\Ykj\right)\\
&&+{\binom {4} {1, 1, 2}}\Esp_{\Phi}\left(\sum_{1\leq i<j<k<l\leq n} \Yij\Yik\Yil+\Yij\Yki\Ykl+\Yij\Yli\Ylk\right) \\
&&+{\binom {5} {2}}\Esp_{\Phi}\left(\sum_{1\leq i<j<k<l<m\leq n} \Yij\Ykl\Ykm+\Yij\Ylk\Ylm+\Yij\Ymk\Ylm\right)  \\
&=& {\binom {n} {2, 1, n-3}}\left(2\phi_2+\phi_3\right)
+{\binom {n} {1, 1, 2, n-4}}\left(\phi_5 + 2\phi_6\right) 
+{\binom {n} {2, 3, n-5}}\left(3\phi_1\phi_2\right)  \\
&=& {\frac{n_2}{2}(2\phi_2+\phi_3) 
+ \frac{n_3}{2}(\phi_5 + 2\phi_6) 
+ \frac{n_4}{4}\phi_1\phi_2}\\
\end{eqnarray*} 
and
\begin{eqnarray*}
\Esp_{\Phi}(M_2^2)
&=&\sum_{1\leq i<j<k\leq n}\Esp_{\Phi}\Big(\Yij^2\Yik^2+\Yji^2\Yjk^2+\Yki^2\Ykj^2\\
&&+2\left(\Yij^2\Yik\Yjk+\Yij\Yik^2\Yjk+\Yij\Yik\Ykj^2\right)\Big)\\
&+&{\binom {4} {2, 1, 1}}\sum_{1\leq i<j<k<l\leq n} \Esp_{\Phi}\Big(\Yij\Yik\Yjk\Yjl+\Yij\Yik\Ykj\Ykl+\Yij\Yik\Ylj\Ylk\\
&&+ \Yji\Yjk^2\Yjl+\Yji\Yjk^2\Ykl+\Yji\Yjk\Ylj\Ylk\\
&&+ \Yki\Ykj^2\Yjl+\Yki\Ykj^2\Ykl+\Yki\Ykj\Ylj\Ylk\Big)\\
&+&{\binom {5} {1, 2, 2}}\sum_{1\leq i<j<k<l<m\leq n} \Esp_{\Phi}\Big(\Yij\Yik\Ykl\Ykm+\Yij\Yik\Ylk\Ylm+\Yij\Yik\Ymk\Ylm\\
&&+\Yji\Yjk\Ykl\Ykm+\Yji\Yjk\Ylk\Ylm+\Yji\Yjk\Ylk\Ylm\\
&&+ \Yki\Ykj\Ykl\Ykm+\Yki\Ykj\Ylk\Ylm+\Yki\Ykj\Ymk\Ylm\Big)\\
&+&{\binom {6} {3, 3}} \sum_{1\leq i<j<k<l<m<u\leq n} 
\Esp_{\Phi}\Big( \Yij\Yik \Ylm\Ylu {+ \Yij\Yik \Ylm\Ymu + \Yij\Yik \Ylu\Ymu} \\ 
& & {+ \Yij\Yjk \Ylm\Ylu + \Yij\Yjk \Ylm\Ymu + \Yij\Yjk \Ylu\Ymu} \\
& & {+ \Yik\Yjk \Ylm\Ylu + \Yik\Yjk \Ylm\Ymu + \Yik\Yjk \Ylu\Ymu} \Big)\\
%
%
&=&\binom {n} {3} (3\phi_2 +  6\phi_3) 
+ \binom {n} {2, 1, 1, n-4} (4\phi_4 + 2\phi_5 + 2\phi_6 + \phi_7) \\
&& + \binom {n} {1, 2, 2, n-5} (4\phi_8 + 4\phi_{10} + \phi_9) 
+ \binom {n} {3, 3, n-6} (9\phi_2^2)  \\
&=&{\frac{n_2}{6}(3\phi_2 +  6\phi_3) 
+ \frac{n_3}{2}(4\phi_4 + 2\phi_5 + 2\phi_6 + \phi_7) 
+ \frac{n_4}{4}(4\phi_8 + 4\phi_{10} + \phi_9) 
+ \frac{n_5}{4}\phi_2^2 }.
\end{eqnarray*} 

\end{document}

%% file: commands.tex
\usepackage{tikz}

\tikzstyle{node}=[inner sep=0]
\tikzstyle{edge}=[-, line width=.5pt]

\newcommand{\Esp}{\mathbb{E}}

\newcommand{\Sd}{\mathbb{S}}
\newcommand{\Var}{\mathbb{V}}
\newcommand{\Bcal}{\mathcal{B}}
\newcommand{\Ecal}{\mathcal{E}}
\newcommand{\Gcal}{\mathcal{G}}

\newcommand{\Ncal}{\mathcal{N}}
\renewcommand{\Pr}{\mathbb{P}}

\newcommand{\Db}{\overline{D}}
\newcommand{\dd}{\text{d}}

\newcommand{\Vh}{V}

\newcommand{\Wo}{W_{\pgo}}
\newcommand{\Who}{W_{\pgoh}}
\newcommand{\Wphio}{W_{\phio}}
\newcommand{\Wphiho}{W_{\widehat{\Phi}^0}}

\newcommand{\Yij}{{Y_{ij}}}
\newcommand{\Yji}{{Y_{ji}}}
\newcommand{\Yik}{{Y_{ik}}}
\newcommand{\Yki}{{Y_{ki}}}
\newcommand{\Yil}{{Y_{i\ell}}}
\newcommand{\Yli}{{Y_{\ell i}}}
\newcommand{\Yjk}{{Y_{jk}}}
\newcommand{\Ykj}{{Y_{kj}}}
\newcommand{\Yjl}{{Y_{j\ell}}}
\newcommand{\Ykl}{{Y_{k\ell}}}
\newcommand{\Ylk}{{Y_{\ell k}}}
\newcommand{\Ylj}{{Y_{\ell j}}}
\newcommand{\Ykm}{{Y_{km}}}

\newcommand{\Ylm}{{Y_{\ell m}}}
\newcommand{\Ylu}{{Y_{\ell u}}}
\newcommand{\Ymk}{{Y_{mk}}}
\newcommand{\Ymu}{{Y_{mu}}}

\newcommand{\Ytij}{{\widetilde{Y}_{ij}}}
\newcommand{\Ytik}{{\widetilde{Y}_{ik}}}

\newcommand{\Ytjk}{{\widetilde{Y}_{jk}}}

\newcommand{\Ytkl}{{\widetilde{Y}_{k\ell}}}

\newcommand{\pb}{{\overline{p}}}
\newcommand{\pg}{{\bf p}}
\newcommand{\pgo}{{\pg^0}}
\newcommand{\pgso}{{\pg^{0*}}}
\newcommand{\pgoh}{{\widehat{\pg}^0}}

\newcommand{\ph}{{\widehat{p}}}
\newcommand{\qh}{{\widehat{q}}}

\newcommand{\phiso}{{\Phi^{0*}}}
\newcommand{\phio}{{\Phi^0}}
\newcommand{\phioo}{{\phi^0_1}}
\newcommand{\phioh}{{\widehat{\Phi}^0}}

\newcommand{\pij}{{p_{ij}}}
\newcommand{\pik}{{p_{ik}}}
\newcommand{\pil}{{p_{i\ell}}}
\newcommand{\pjk}{{p_{jk}}}
\newcommand{\pjl}{{p_{j\ell}}}
\newcommand{\pkl}{{p_{k\ell}}}

\newcommand{\poij}{{p^0_{ij}}}

\newcommand{\phoij}{{\widehat{p}^0_{ij}}}

\newcommand{\psij}{{p^*_{ij}}}
\newcommand{\psoij}{{p^{0*}_{ij}}}

\newcommand{\Di}{{\Delta_i}}
\newcommand{\Dj}{{\Delta_j}}

\newcommand{\dij}{{\delta_{ij}}}
\newcommand{\dik}{{\delta_{ik}}}

\newcommand{\djk}{{\delta_{jk}}}

\newcommand{\dcij}{{\delta^2_{ij}}}

\newcommand{\sumij}{{\sum_{i \neq j}}}
\newcommand{\sumji}{{\sum_{j \neq i}}}

\newcommand{\sigmaij}{{\sigma^2_{ij}}}
\newcommand{\sigmaik}{{\sigma^2_{ik}}}
\newcommand{\sigmail}{{\sigma^2_{i\ell}}}
\newcommand{\sigmajk}{{\sigma^2_{jk}}}
\newcommand{\sigmajl}{{\sigma^2_{j\ell}}}
\newcommand{\sigmakl}{{\sigma^2_{k\ell}}}

\newcommand{\sigmaoij}{{{\sigma^0_{ij}}^2}}
\newcommand{\sigmaoik}{{{\sigma^0_{ik}}^2}}

\newcommand{\sigmaojk}{{{\sigma^0_{jk}}^2}}

\newcommand{\xij}{{x_{ij}}}
\newcommand{\xpij}{{x^\intercal_{ij}}}
\newcommand{\xoij}{{x^0_{ij}}}
\newcommand{\xopij}{{x^{0\intercal}_{ij}}}

\newcommand{\anu}{a_{nu}}
\newcommand{\pnu}{p_{nu}}
\newcommand{\Xnu}{X_{nu}}
\newcommand{\Znu}{Z_{nu}}
\newcommand{\sigmanu}{\sigma_{nu}}

\newcommand{\ERph}{{ER(\ph)}}

\newcommand{\HERpg}{{HER(\pg)}}
\newcommand{\HERpgo}{{HER(\pgo)}}

\newcommand{\EG}{{EG}}
\newcommand{\EGphi}{{\EG(\Phi)}}
\newcommand{\EGphio}{{\EG(\Phi^0)}}
\newcommand{\EGphioh}{{\EG(\widehat{\Phi}^0)}}

\newcommand{\muo}{{\mu}^0}
\newcommand{\muho}{\widehat{\mu}^0}

\newtheorem{assumption}{Assumption}
\newtheorem{corollary}{Corollary}

\newtheorem{lemma}{Lemma}
\newtheorem{proposition}{Proposition}
\newtheorem{remark}{Remark}
\newtheorem{theorem}{Theorem}

\newcommand{\proofbegin}{\noindent{\sl Proof.}~}
\newcommand{\proofend}{$\blacksquare$\bigskip}

%% file: article_arXiv.bbl
\begin{thebibliography}{}

\bibitem[\protect\citeauthoryear{Barab\'asi and Albert}{Barab\'asi and
  Albert}{1999}]{BaA99}
Barab\'asi, A.~L. and R.~Albert (1999).
\newblock Emergence of scaling in random networks.
\newblock {\em Science\/}~{\em 286}, 509--512.

\bibitem[\protect\citeauthoryear{Bickel, Chen, and Levina}{Bickel
  et~al.}{2011}]{BCL11}
Bickel, P.~J., A.~Chen, and E.~Levina (2011).
\newblock The method of moments and degree distributions for network models.
\newblock {\em Ann. Stat.\/}~{\em 39\/}(5), 2280--2301.

\bibitem[\protect\citeauthoryear{Bickel and Sarkar}{Bickel and
  Sarkar}{2016}]{BiSa16}
Bickel, P.~J. and P.~Sarkar (2016).
\newblock Hypothesis testing for automated community detection in networks.
\newblock {\em Journal of the Royal Statistical Society: Series B (Statistical
  Methodology)\/}~{\em 78\/}(1), 253--273.

\bibitem[\protect\citeauthoryear{Billingsley}{Billingsley}{1968}]{Bil68}
Billingsley, P. (1968).
\newblock {\em Convergence of Probability Measures}.
\newblock Wiley: New-York.

\bibitem[\protect\citeauthoryear{Bloznelis}{Bloznelis}{2005}]{Blo05}
Bloznelis, M. (2005).
\newblock Degree variance is asymptotically normal.
\newblock Technical report, Vilnius university, Faculty of Mathematics and
  Informatics.

\bibitem[\protect\citeauthoryear{Burt}{Burt}{1987}]{Burt1987}
Burt, R. (1987).
\newblock Social contagion and innovation: cohesion versus structural
  equivalence.
\newblock {\em American Journal of Sociology\/}~{\em 92}, 1287--1335.

\bibitem[\protect\citeauthoryear{Cerqueira, Fraiman, Vargas, and
  Leonardi}{Cerqueira et~al.}{2017}]{CFV15}
Cerqueira, A., D.~Fraiman, C.~D. Vargas, and F.~Leonardi (2017).
\newblock A test of hypotheses for random graph distributions built from eeg
  data.
\newblock {\em IEEE Transactions on Network Science and Engineering\/}~{\em
  4\/}(2), 75--82.

\bibitem[\protect\citeauthoryear{{Channarond}, {Daudin}, and
  {Robin}}{{Channarond} et~al.}{2012}]{CDR12}
{Channarond}, A., J.-J. {Daudin}, and S.~{Robin} (2012).
\newblock Classification and estimation in the stochastic block model based on
  the empirical degrees.
\newblock {\em Elec. J. Stat.\/}~{\em 6}, 2574--601.

\bibitem[\protect\citeauthoryear{Chung and Lu}{Chung and Lu}{2002}]{ChL02}
Chung, F. and L.~Lu (2002).
\newblock Connected components in random graphs with given expected degree
  sequences.
\newblock {\em Annals of combinatorics\/}~{\em 6\/}(2), 125--145.

\bibitem[\protect\citeauthoryear{Coleman, Katz, and Menzel}{Coleman
  et~al.}{1966}]{Coleman1966}
Coleman, J., E.~Katz, and H.~Menzel (1966).
\newblock Medical innovation: a diffusion study. indianapolis: the
  boobs-merrill company.
\newblock {\em Behavioral Science\/}~{\em 12}, 481--483.

\bibitem[\protect\citeauthoryear{Dasgupta, Hopcroft, and McSherry}{Dasgupta
  et~al.}{2004}]{dasgupta2004spectral}
Dasgupta, A., J.~E. Hopcroft, and F.~McSherry (2004).
\newblock Spectral analysis of random graphs with skewed degree distributions.
\newblock In {\em null}, pp.\  602--610. IEEE.

\bibitem[\protect\citeauthoryear{Daudin, Picard, and Robin}{Daudin
  et~al.}{2008}]{DPR08}
Daudin, J.-J., F.~Picard, and S.~Robin (2008).
\newblock A mixture model for random graphs.
\newblock {\em Stat. Comput.\/}~{\em 18\/}(2), 173--83.

\bibitem[\protect\citeauthoryear{Diaconis and Janson}{Diaconis and
  Janson}{2008}]{DiJ08}
Diaconis, P. and S.~Janson (2008).
\newblock Graph limits and exchangeable random graphs.
\newblock {\em {R}end. {M}at. {A}ppl.\/}~{\em 7\/}(28), 33--61.

\bibitem[\protect\citeauthoryear{Erd\"os and R\'enyi}{Erd\"os and
  R\'enyi}{1959}]{ER59}
Erd\"os, P. and A.~R\'enyi (1959).
\newblock On random graphs.
\newblock {\em I Publicationes Mathematicae (Debrecen)\/}~{\em 6}, 290--297.

\bibitem[\protect\citeauthoryear{Gao and Lafferty}{Gao and
  Lafferty}{2017a}]{GaL17b}
Gao, C. and J.~Lafferty (2017a).
\newblock Testing for global network structure using small subgraph statistics.
\newblock Technical Report 1710.00862, arXiv.

\bibitem[\protect\citeauthoryear{Gao and Lafferty}{Gao and
  Lafferty}{2017b}]{Gao17}
Gao, C. and J.~Lafferty (2017b).
\newblock Testing network structure using relations between small subgraph
  probabilities.
\newblock {\em arXiv preprint arXiv:1704.06742\/}.

\bibitem[\protect\citeauthoryear{Hagberg}{Hagberg}{2000}]{Hag00}
Hagberg, J. (2000).
\newblock Centrality testing and the distribution of the degree variance in
  bernoulli graphs.
\newblock Technical report, Department of Statistics, Stockholm University.

\bibitem[\protect\citeauthoryear{Hagberg}{Hagberg}{2003}]{Hag03a}
Hagberg, J. (2003).
\newblock General moments of degrees in random graphs.
\newblock {\em Stockholm University, Department of Statistics\/}.

\bibitem[\protect\citeauthoryear{Handcock, Hunter, Butss, Goodreau, and
  Morris}{Handcock et~al.}{2008}]{handcock2008}
Handcock, M., D.~Hunter, C.~Butss, S.~Goodreau, and M.~Morris (2008).
\newblock Statnet: Software tools for the representation, visualization,
  analysis and simulation of network data.
\newblock {\em Journal of Statistical Software\/}~{\em 24}, 12--25.

\bibitem[\protect\citeauthoryear{Holland and Leinhardt}{Holland and
  Leinhardt}{1979}]{Hol79}
Holland, P.~W. and S.~Leinhardt (1979).
\newblock Structural sociometry.
\newblock {\em Perspectives on social network research\/}, 63--83.

\bibitem[\protect\citeauthoryear{Hunter, Goodreau, and Handcock}{Hunter
  et~al.}{2008}]{HGH08}
Hunter, D.~R., S.~M. Goodreau, and M.~S. Handcock (2008).
\newblock Goodness of fit of social network models.
\newblock {\em Journal of the American Statistical Association\/}~{\em
  103\/}(481), 248--258.

\bibitem[\protect\citeauthoryear{Josse and Husson}{Josse and
  Husson}{2016}]{husson2016}
Josse, J. and F.~Husson (2016).
\newblock {missMDA}: a package for handling missing values in multivariate data
  analysis.
\newblock {\em Journal of Statistical Software\/}~{\em 70\/}(1), 1--31.

\bibitem[\protect\citeauthoryear{Karrer and Newman}{Karrer and
  Newman}{2011}]{karrer11}
Karrer, B. and M.~E.~J. Newman (2011).
\newblock Stochastic blockmodels and community structure in networks.
\newblock {\em Phys. Rev. E\/}~{\em 83}, 016107.

\bibitem[\protect\citeauthoryear{Latouche, Birmel\'e, and Ambroise}{Latouche
  et~al.}{2011}]{LBA11a}
Latouche, P., E.~Birmel\'e, and C.~Ambroise (2011).
\newblock {Overlapping stochastic block models with application to the French
  political blogosphere.}
\newblock {\em Ann. Appl. Stat.\/}~{\em 5\/}(1), 309--336.

\bibitem[\protect\citeauthoryear{Latouche and Robin}{Latouche and
  Robin}{2016}]{LaR15}
Latouche, P. and S.~Robin (2016).
\newblock Variational bayes model averaging for graphon functions and motif
  frequencies inference in {$W$}-graph models.
\newblock {\em Statistics and Computing\/}~{\em 26}, 1173--1185.

\bibitem[\protect\citeauthoryear{Latouche, Robin, and Ouadah}{Latouche
  et~al.}{2018}]{gof2}
Latouche, P., S.~Robin, and S.~Ouadah (2018).
\newblock Goodness of fit of logistic models for random graphs.
\newblock {\em Journal of Computational and Graphical Statistics\/}~{\em
  27\/}(1), 98--109.

\bibitem[\protect\citeauthoryear{Lei}{Lei}{2016}]{Lei16}
Lei, J. (2016).
\newblock A goodness-of-fit test for stochastic block models.
\newblock {\em The Annals of Statistics\/}~{\em 44\/}(1), 401--424.

\bibitem[\protect\citeauthoryear{Lov\'asz and Szegedy}{Lov\'asz and
  Szegedy}{2006}]{LoS06}
Lov\'asz, L. and B.~Szegedy (2006).
\newblock Limits of dense graph sequences.
\newblock {\em Journal of Combinatorial Theory, Series B\/}~{\em 96\/}(6), 933
  -- 957.

\bibitem[\protect\citeauthoryear{Mariadassou, Robin, and Vacher}{Mariadassou
  et~al.}{2010}]{MRV10}
Mariadassou, M., S.~Robin, and C.~Vacher (2010).
\newblock Uncovering structure in valued graphs: a variational approach.
\newblock {\em Ann. Appl. Statist.\/}~{\em 4\/}(2), 715--42.

\bibitem[\protect\citeauthoryear{Maugis, Priebe, Olhede, and Wolfe}{Maugis
  et~al.}{2017}]{Mau17}
Maugis, P., C.~E. Priebe, S.~C. Olhede, and P.~J. Wolfe (2017).
\newblock Statistical inference for network samples using subgraph counts.
\newblock {\em arXiv preprint arXiv:1701.00505\/}.

\bibitem[\protect\citeauthoryear{Newman}{Newman}{2003}]{New03}
Newman, M.~E. (2003).
\newblock The structure and function of complex networks.
\newblock {\em SIAM review\/}~{\em 45\/}(2), 167--256.

\bibitem[\protect\citeauthoryear{Nowicki and Snijders}{Nowicki and
  Snijders}{2001}]{NoS01}
Nowicki, K. and T.~Snijders (2001).
\newblock Estimation and prediction for stochastic block-structures.
\newblock {\em J. Amer. Statist. Ass.\/}~{\em 96}, 1077--87.

\bibitem[\protect\citeauthoryear{Nowicki and Wierman}{Nowicki and
  Wierman}{1988}]{NoW88}
Nowicki, K. and J.~C. Wierman (1988).
\newblock Subgraph counts in random graphs using incomplete u-statistics
  methods.
\newblock {\em Discrete Math.\/}~{\em 72\/}(1), 299--310.

\bibitem[\protect\citeauthoryear{Picard, Daudin, Koskas, Schbath, and
  Robin}{Picard et~al.}{2008}]{PDK08}
Picard, F., J.-J. Daudin, M.~Koskas, S.~Schbath, and S.~Robin (2008).
\newblock Assessing the exceptionality of network motifs,.
\newblock {\em J. Comput. Biol.\/}~{\em 15\/}(1), 1--20.

\bibitem[\protect\citeauthoryear{Rasch}{Rasch}{1960}]{rasch1960probabilistic}
Rasch, G. (1960).
\newblock {\em Probabilistic Models for Some Intelligence and Attainment
  Tests}.
\newblock Studies in mathematical psychology. Danmarks Paedagogiske Institut.

\bibitem[\protect\citeauthoryear{Resnick, Bearman, Blum, Bauman, Harris, Jones,
  Tabor, Beuhring, Sieving, Shew, et~al.}{Resnick et~al.}{1997}]{Resnick97}
Resnick, M., P.~S. Bearman, R.~W. Blum, K.~E. Bauman, K.~M. Harris, J.~Jones,
  J.~Tabor, T.~Beuhring, R.~E. Sieving, M.~Shew, et~al. (1997).
\newblock Protecting adolescents from harm: findings from the national
  longitudinal study on adolescent health.
\newblock {\em Jama\/}~{\em 278\/}(10), 823--832.

\bibitem[\protect\citeauthoryear{Snijders}{Snijders}{1981}]{Sni81}
Snijders, T. A.~B. (1981).
\newblock The degree variance: An index of graph heterogeneity.
\newblock {\em Social Networks\/}~{\em 3\/}(3), 163--174.

\bibitem[\protect\citeauthoryear{Vacher, Piou, and Desprez-Loustau}{Vacher
  et~al.}{2008}]{VPD08}
Vacher, C., D.~Piou, and M.-L. Desprez-Loustau (2008).
\newblock Architecture of an antagonistic tree/fungus network: The asymmetric
  influence of past evolutionary history.
\newblock {\em PLoS ONE\/}~{\em 3\/}(3), 1740.

\bibitem[\protect\citeauthoryear{van~der Vaart}{van~der Vaart}{1998}]{vdV98}
van~der Vaart, A.~W. (1998).
\newblock {\em Asymptotic statistics}, Volume~3 of {\em Cambridge Series in
  Statistical and Probabilistic Mathematics}.
\newblock Cambridge University Press, Cambridge.

\bibitem[\protect\citeauthoryear{Yang, Han, and Airoldi}{Yang
  et~al.}{2014}]{YHA14}
Yang, J., C.~Han, and E.~Airoldi (2014).
\newblock Nonparametric estimation and testing of exchangeable graph models.
\newblock In {\em AISTATS}, pp.\  1060--1067.

\bibitem[\protect\citeauthoryear{Young and Scheinerman}{Young and
  Scheinerman}{2007}]{YoS07}
Young, S.~J. and E.~R. Scheinerman (2007).
\newblock Random dot product graph models for social networks.
\newblock In {\em International Workshop on Algorithms and Models for the
  Web-Graph}, pp.\  138--149. Springer.

\bibitem[\protect\citeauthoryear{Zachary}{Zachary}{1977}]{zachary1977}
Zachary, W. (1977).
\newblock An information flow model for conflict and fission in small groups.
\newblock {\em Journal of Anthropological Research\/}~{\em 33}, 452--473.

\end{thebibliography}
